\documentclass[twoside,12pt]{article}
\usepackage[abbrvbib, preprint]{jmlr2e}
\setcitestyle{square,comma, open={[}, close={]}}
\setcitestyle{numbers}
\usepackage{hyperref}
\usepackage{tabularx}
\usepackage{amsmath, amssymb}
\usepackage{hyperref}
\hypersetup{hidelinks}
\usepackage[ruled,vlined]{algorithm2e}
\usepackage{algorithmic}
\usepackage{mathrsfs}
\usepackage{pgfplots}
\usepackage[skip=2pt,font=scriptsize]{caption}
\usepackage{float}
\usepackage{subfig} 
\usepackage{bbm}
\usepackage{bm}
\usepackage[normalem]{ulem}

\usepackage{lipsum}
\usepackage [english]{babel}
\usepackage [autostyle, english = american]{csquotes}
\usepackage{soul}
\usepackage{graphicx}
\usepackage{csvsimple}
\usepackage{booktabs}  
\usepackage{pifont}  
\usepackage{xcolor} 
\usepackage{enumitem} 

\let\proof\relax

\usepackage{amsthm, mathtools, thmtools}

\usepackage{cleveref}

{\bfseries}{\itshape}
{\bfseries}{\itshape}
{\bfseries}{\itshape}
{\bfseries}{\itshape}
{\bfseries}{\itshape}
{\bfseries}{\itshape}
{\bfseries}{\itshape}

\DeclareMathOperator{\E}{\mathbb{E}}

\DeclareMathOperator{\Var}{Var}

\DeclareMathOperator{\Span}{span}
\DeclareMathOperator{\Aff}{aff}
\DeclareMathOperator{\Int}{int}
\DeclareMathOperator{\ri}{ri}
\DeclareMathOperator{\cm}{cm}
\DeclareMathOperator{\bd}{bd}

\DeclareMathOperator{\conv}{conv}

\newcommand{\ip}[2]{\left\langle #1,#2\right\rangle}

\newcommand{\ppphi}{\mathcal{A}}

\newcommand{\R}{\mathbb{R}}
\newcommand{\Pp}{\mathbb{P}}
\newcommand{\norm}[1]{\left\|#1\right\|_2}

\def \C {\mathbb{C}}

\def \P {\mathbb{P}}
\def \R {\mathbb{R}}

\def \RR {\mathcal{R}}

\renewcommand{\epsilon}{\varepsilon}

\newtheorem{theorem}{Theorem}[section]

\newtheorem{corollary}[theorem]{Corollary}
\newtheorem{lemma}[theorem]{Lemma}

\newtheorem{definition}[theorem]{Definition}

\definecolor{darkgreen}{rgb}{0.3, 0.7, 0.3}

\def\RR{\mathbb{R}}
\def\ll{\langle}

\def\rr{\rangle}

\begin{document}

\title{\bf Conic Frames, Cosine Measure, and Gaussian Frames }

\author{\name Jonathan Forstater \email jforstater@ucdavis.edu \\
\addr Department of Mathematics \\ University of California, Davis \\ 
       Davis, CA 95616-5270, USA
       \AND
       \name Thomas Strohmer \email strohmer@math.ucdavis.edu \\
       \addr Department of Mathematics\\
       University of California Davis\\
       Davis, CA 95616-5270, USA}

\maketitle

\begin{abstract}
Conic combinations play an important role across mathematics in areas such as optimization, non-negative matrix factorization, and sparse coding. Frames which span their underlying space using only conic combinations, known as positive spanning sets or conic frames, are of special interest in derivative free optimization. In particular, sets with large cosine measure, a measure of how spread out the vectors are, are favored. We focus on the asymptotics of the cosine measure of Gaussian frames and present tail bounds in a variety of regimes.
\end{abstract}

\section{Introduction}

Frames are a fundamental concept in mathematics, particularly in areas such as harmonic analysis, linear algebra, functional analysis, signal processing, and numerical analysis~\cite{young1981introduction,daubechies1992ten,christensen2008frames}. They generalize the idea of a basis in a vector space, offering a more flexible and robust way to represent and analyze signals or data.
Frames offer a powerful mathematical framework for understanding and manipulating data in high-dimensional spaces, providing a bridge between the rigid structure of bases and the flexibility needed for real-world applications. As such, they have found widespread applications in signal and image processing (e.g., wavelet frames for multiresolution analysis, Gabor frames for time-frequency analysis), data compression, machine learning, and numerical methods for solving partial differential equations. 

Unlike a basis, a frame can be redundant, meaning its vectors are not necessarily linearly independent.
Indeed, this non-uniqueness of representation coefficients turns out to be very useful for certain applications where it is desirable to restrict the coefficients to be {\em non-negative}, resulting in what is known as a {\em conic combination}. This type of constraint is fundamental to the mathematical field of optimization and has applications in many scientific and engineering problems, such as Non-Negative Matrix Factorization, Sparse Coding and Dictionary Learning, and Derivative Free Optimization.

We organize the paper as follows. In \Cref{sec:background}, we review the necessary background on frames, convex geometry, and their applications. Then, in \Cref{sec:existence-construction}, we develop geometric characterizations of conic frames, establish their connection with bounded polyhedra, and present randomized and deterministic methods to make certain frames conic. In \Cref{sec:max-cosine-measure}, we examine the relationship between coherence and the cosine measure and show that twin equiangular tight frames and regular polytopes need not maximize the cosine measure. Last, in \Cref{sec:cm-gaussianframes}, we identify the cosine measure with the minimum facet height of the convex hull of the frame and use this connection to establish semi-asymptotic formulas and concentration bounds for normalized Gaussian frames. We conclude by providing non-asymptotic tail bounds in the subexponential regime.

\section{Mathematical Background}
\label{sec:background}

\subsection{Frames}

We  recommend~\cite{christensen2008frames} as a great introduction to frame theory.
A frame in a separable Hilbert space $\mathcal{H}$  is a sequence of vectors $\{f_k\}_{k \in I}$ (where $I$ is an index set, possibly infinite) such that there exist constants $A > 0$ and $B < \infty$ satisfying the frame bounds for all vectors $f \in \mathcal{H}$,
\begin{equation}\label{framebound}
A\|f\|^2 \leq \sum_{k \in I} |\langle f, f_k \rangle|^2 \leq B\|f\|^2.
\end{equation}

Each frame has the associated frame operator, $S: \mathcal{H} \to \mathcal{H}$:
$$Sf = \sum\limits_{k \in I} \langle f, f_k \rangle f_k.$$
Note that the frame operator is self-adjoint and the lower frame bound implies that $S$ is (strictly) positive and thus invertible. This inverse allows us to define the canonical dual frame, $\{S^{-1}f_k\}_{k \in I}$, leading to the following two representations:

\begin{align*}
    f = \sum\limits_{k \in I} \langle f, f_k \rangle S^{-1}f_k = \sum\limits_{k \in I} \langle f, S^{-1}f_k \rangle f_k
\end{align*}
If the frame is tight, $S$ and $S^{-1}$ are both multiples of the identity, meaning that computing these representations are efficient.

In this paper we focus on {\em finite frames}, i.e., frames for a finite-dimensional Hilbert space~\cite{christensen2008frames,waldron2018introduction,casazza2012finite}. In particular, we are interested in the case ${\mathcal H} = \R^n$. Thus, we consider vectors $f_1,\dots,f_m$ with $m \ge n$. The upper frame bound is trivially satisfied in this case. To satisfy the lower frame bound in~\eqref{framebound} the vectors $f_1,\dots,f_m$ need to span $\R^n$. Equivalently, we can require that the $n \times m$ matrix $[f_1,\dots,f_m]$ has rank $n$. By a slight abuse of notation, we also write $F$ for the matrix representing a frame, i.e.,
$F=[f_1,\dots,f_m]$.

For certain applications, we are interested in how `spread out' the frame vectors are. This is quantified by the notion of \emph{coherence}:

\begin{definition}[Coherence]
    The coherence, $\mu$, of a frame $F$ consisting of all non-zero vectors is: $$\mu(F) = \max\limits_{k \neq l} \frac{| \langle f_k, f_l \rangle |}{\| f_k \| \| f_l \|}$$
\end{definition}

Notice that coherence is equivalent to the minimum acute angle between the lines spanned by the frame vectors. By minimizing the coherence, the `spread' of the vectors is thus maximized. The following lower bound on the coherence was originally shown by L. R. Welch~\cite{welch_lower_1974}:

\begin{theorem}[Welch Bound]
    Let $\{ f_k \}_{k = 1}^m$ be a unit norm frame for $\R^n$ or $\C^n$, then $$\mu \ge \sqrt{\frac{m - n}{n(m - 1)}}$$
\end{theorem}

Tightness of a frame is a necessary condition to achieve the Welch bound, but it is not sufficient; this leads us to our next definition

\begin{definition}[Equiangularity]
    A frame is called equiangular if $| \langle f_k, f_l \rangle| = \alpha$ for some fixed $\alpha$ and for all $k \neq l$.
\end{definition}

This leads us to  {\em equiangular tight frames} (ETF), which rose to prominence via~\cite{StrohmerHeath}. It is proven there that equality holds in the Welch bound if and only if a frame is an ETF.

\subsection{Conic Frames and Convex Geometry}

Conic frames are frames that can represent every vector in their underlying space as a conic combination (a linear combination with non-negative coefficients). 

\begin{definition}[Conic Frame]
    Let $F \in \R^{m \times n}$ be a finite frame such that for all $y \in \R^n$ there exists an $x \in \R^m$ with $x \ge 0$ such that $Fx = y$. Then $F$ is a \emph{conic frame}.
\end{definition}

These are fundamental to convex geometry, along with affine (coefficients sum to one) and convex combinations (both affine and conic).

A natural way to study frames 
is by studying the properties of their convex hulls, which are known as convex polytopes. In order to do so, we will first provide some background definitions and results from convex geometry. 

Let $M \subseteq \R^n$. The affine hull (or affine span) of $M$, $\Aff M$, is the smallest affine set that contains $M$. Similarly, the convex hull of $M$, $\conv M$ is the smallest convex set containing $M$.

Affine sets are in fact translates of linear subspaces (also known as a \emph{flats}). That is, if $A$ is an affine set then $A = x + L$, where $x \in \R^n$ and $L$ is a linear subspace. The dimension of an affine set is defined as the dimension of $L$, while the dimension of a convex set, $C$, is the dimension of $\Aff C$.

The following lemma (the proof of which is left to the reader) relates the affine span of a set to its linear span. 

\begin{lemma}\label{lemma:dim-aff}
    For any subset $M$ of $\R^n$, if $0 \in \Aff M$, then $\dim(\Aff M) = \dim(\Span M)$. Otherwise, $\dim(\Aff M) = \dim(\Span M) - 1$
\end{lemma}

Convex sets that have a lower dimension than their ambient space have an empty interior. This leads us to the notion of \emph{relative interior} which is guaranteed to exist for any non-empty convex set:

\begin{definition}
    The \emph{relative interior} of a convex set $C$, $\ri C$ is the interior of $C$ in $\Aff C$.
\end{definition}

Note that a non-empty convex set in $C \subseteq \R^n$ has a non-empty interior if and only if $\Aff C = \R^n$ and in such a case $\ri C = \Int C$~\cite{brondsted_introduction_1983}.

As mentioned previously, a type of convex set that is particularly useful in understanding finite frames is the \emph{convex polytope}, which is the convex hull of a finite set of points. Note that a convex polytope is closed and bounded and is thus a compact convex set. From here on, we may refer to convex polytopes simply as polytopes. The following lemma (which the reader can readily verify) allows us to characterize the relative interior of a polytope. 

\begin{lemma}\label{lemma:ri-iff-strictpos}
    Let $P = \conv\{x_1, \dots, x_m\}$ be a polytope in $\R^n$. $x \in \ri P$ if and only if $x$ is a convex combination of $\{x_1, \dots, x_m\}$ with strictly positive coefficients.
\end{lemma}

In order to understand the structure of a polytope, we can look at its \emph{faces}. A \emph{supporting hyperplane} of a convex set is a hyperplane which contains the convex set in one of its associated halfspaces and intersects the set; a face is precisely this intersection.

An important feature of the faces of polytopes is that they themselves are polytopes. The following is a rephrasing of Proposition 2.3 in~\cite{ziegler_lectures_1995}:

\begin{theorem}[Face of Convex Polytope is Convex Polytope]\label{thm:face-poly}
    Let $P$ be a convex polytope, then each face of $P$ is itself a convex polytope. Specifically, if $A = \{ x_1, ...., x_k \}$ and $P = conv (A)$, a face of $P$ has the form $\conv (B)$ where $B \subseteq A$.
\end{theorem}

Noting that a hyperplane is $(n-1)$-dimensional and every face is a subset of some hyperplane, we have the following corollary:

\begin{corollary}[Faces of Polytopes are at most $n-1$-dimensional]\label{cor:faces-n-1}
     Let $K$ be an $n$-dimensional polytope embedded in $\RR^n$, then each face of $K$ is a $k$-dimensional polytope (with $k < n$), defined by a subset of the points that define $K$.
\end{corollary}

The following two theorems give conditions for the existence of supporting hyperplanes. They can be found in Sections 2.5.1 and 2.5.2 of~\cite{boyd2004convex} respectively where they are phrased slightly differently.

\begin{theorem}[Separating Hyperplane Theorem]\label{thm:sep-hyperplane}
    Let $C$ and $D$ be non-empty disjoint convex sets in $\R^n$. Then there is a hyperplane $H$ such that $C \subseteq H^+$ and $D \subseteq H^-$.
\end{theorem}

\begin{theorem}[Supporting Hyperplane Theorem]
    Let $C$ be a non-empty closed convex set in $\R^n$ and let $b$ be a boundary point of $C$. Then there exists at least one supporting hyperplane of $C$ containing $b$.
\end{theorem}

This last theorem can be restated in a way that will be useful later in our study of conic frames:

\begin{corollary}\label{cor:boundary-face}
    Each boundary point of a closed convex set $C$ in $\RR^n$ lies in at least one face of $C$.
\end{corollary}

\subsection{Previous Work on Conic Combinations and their Applications}

Chandler Davis' paper~\cite{davis1954theory} is one of the starting points for studying conic combinations\footnote{The term {\em frame} also appears in the paper by Chandler Davis, however, it has a completely different meaning there and should not be confused with the widely established concept of a frame as we use it in this paper}. In that paper, he refers to conic combinations as positive combinations and defines notions of positive independence, positive span, and positive basis, which mirror their linear versions. He then mainly focuses on discussing positive bases, which are guaranteed to have between $n+1$ and $2n$ elements in $\R^n$. Follow-up work by~\cite{mckinney1962positive} generalizes some results to infinite dimensional spaces and develops the theory of positive bases further.

At this point, the term \emph{positive spanning set} came into use to refer to the object we now refer to as a conic frame. In \cite{regis2016properties}, Regis discusses positive spanning sets, positive independent sets, and positive bases further. It also shows how many of the things we take for granted about linear versions of such sets are not necessarily true for the positive versions. He also discusses algorithms for extending any set of vectors in a linear subspace to a positive spanning set of that subspace by adding additional vectors. In \Cref{sec:existence-construction}, we present deterministic and randomized algorithms to turn a full-spark frame into a conic frame by merely changing the signs on the vectors of the frame without adding vectors.

Conic combinations appear in a wide range of scientific and engineering disciplines due to their role in describing positive linear systems, optimization problems, and geometric structures. In particular, conic combinations underpin the theory of convex cones, which are central to many areas of applied mathematics. A fundamental example is found in linear programming and its generalization, conic programming, where feasible regions are often defined using conic constraints—sets closed under conic combinations~\cite{boyd2004convex}.

Non-negative Matrix Factorization (NMF) is a prime example where conic combinations are implicitly central~\cite{daniel1999learning}. NMF decomposes a non-negative data matrix into two non-negative matrices. The columns of the data matrix are approximated as conic combinations of the columns of one of the factor matrices (the "basis vectors" or "features"). 

In Audio Source Separation~\cite{virtanen_compositional_2015}, mixed audio signals are decomposed into individual sound sources using non-negative combinations of atoms with NMF techiniques~\cite{virtanen_monaural_2007} or probability distributions in a technique known as probabilistic latent component analysis (PLCA)~\cite{smaragdis2006probabilistic}.

In Beamforming~\cite{van1988beamforming}, an array of antennas is arranged so that by changing the phase and relative amplitude of the signal sent from each antenna, the signal can be formed into a directed beam. Thus, the signal can be viewed as coefficients in a linear system. In general, the coefficients can be complex-valued, but by reformulating this system in terms of intensity, the parameters become non-negative real numbers, leading to a conic combination~\cite{beamforming_conic}.

Conic combinations also play an important role in finance in portfolio theory. When constructing investment portfolios, one often considers non-negative weights for different assets. Imposing such a condition can reduce risk in estimated optimal portfolios~\cite{jagannathan_risk_2003}. 

Another domain where conic combinations are relevant is machine learning, particularly in areas such as sparse coding / dictionary learning~\cite{aharon_k-svd_2006} and support vector machines (SVMs). In sparse coding, an over-complete dictionary of atoms is learned which can be used to represent signals by means of a sparse linear combination of these atoms. Hoyer~\cite{hoyer_non-negative_2002} approximates nonnegative data using learned nonnegative dictionary atoms and nonnegative coefficients, incorporating a sparsity penalty to encourage representations with few active atoms. As for SVMs, Potluru et al.~\cite{potluru_sparseness_2011}
reduce totally nonnegative least squares to a single-class
SVM, relating the support of the nonnegative coefficient
vector to the SVM support vectors.

In economics and game theory, conic combinations help model production sets or strategy spaces where only nonnegative actions or resources are meaningful. This modeling is especially prevalent in input-output analysis and cost modeling, where outputs depend linearly on inputs with nonnegative coefficients~\cite{samuelson1948foundations}.

But perhaps one of the most natural applications of conic combinations and conic frames in particular is Derivative-free Optimization (DFO). In many real-life situations, gradient information is unavailable or costly to compute. In these cases a DFO method, such as directional direct 
search (DDS) is particularly useful. In DDS~\cite{LMW2019AN}, a descent direction is found by polling objective values in a number of directions. These directions correspond to a positive spanning set, or conic frame as it is called here. Under mild assumptions, directional direct search is guaranteed to converge~\cite{conn_introduction_2009}. 
In proving convergence guarantees, the cosine measure of a conic frame is of particular importance.

\begin{definition}[Cosine Measure]
    The cosine measure of a set of nonzero vectors, F, is defined as

    $$\cm F = \min_{0 \neq v \in \R^n} \max_{f \in F} \frac{v^\top f}{\|v\| \|f\|}$$
\end{definition}

The cosine measure of a PSS is guaranteed to fall between 0 and 1 and PSSs with maximal cosine measure are often desired for convergence purposes. In \cite{naevdal2019positive}, N\ae vdal finds positive bases with $n+1$ and $2n$ elements with maximal cosine measure, but for the vast majority of conic frames in $\R^n$ the question of which set of $m$ elements has maximal cosine measure remains open~\cite{dodangeh_optimal_2016}.

\section{Existence and construction of conic frames}
\label{sec:existence-construction}

With the above definitions and results in hand, we are ready to prove some results about conic frames. We begin with necessary and sufficient conditions for a frame to be a conic frame and then connect conic frames with bounded polyhedra. Finally, we present methods to construct a conic frame \emph{de novo} or from an existing full-spark frame with high probability or deterministically.

\subsection{Conic Frames and the Origin}
Some of the results in this section have been shown previously by Regis~\cite{regis2016properties} and Davis~\cite{davis1954theory}, but we present new proofs using polytopes to provide geometric intuition. We begin by demonstrating that the interior of the convex hull of a frame containing the origin is equivalent to certain conditions.

\begin{lemma}\label{lemma:origin-strict-pos}
    Let $F = \{x_1, \dots, x_m\}$ be a frame in $\R^n$ and let $P = \conv(F)$. $0 \in int(P)$ if and only if 0 can be written as a strictly positive convex combination of $x_1, \dots, x_m$.
\end{lemma}

\begin{proof} 
 $(\implies)$ Assume $0 \in \Int P$. Then $\Int P \neq \emptyset$, so $0 \in \ri P = \Int P$ and, by \Cref{lemma:ri-iff-strictpos}, 0 can be written by a strictly positive convex combination of the frame vectors. 
 
 $(\impliedby)$ Let 0 be representable as a strictly positive convex combination of the frame vectors. Then by \Cref{lemma:ri-iff-strictpos}, $0 \in \ri P$. It suffices to show that $\Aff P = \R^n$ or equivalently $\dim \Aff P = n$. Note $0 \in \ri P \subseteq \Aff P$. Thus by \Cref{lemma:dim-aff} we have $\dim(\Aff P) = \dim(\Span P) = n$, where the last equality follows since finite frames are spanning sets.
\end{proof}

\begin{lemma}\label{lemma:origin-hyperplane}
    Let $F$ be a frame in $\R^n$ and $P = \conv(F)$. Then $0 \notin \Int P$ if and only if there is a hyperplane passing through $0$ such that $F$ is contained in one of the resulting halfspaces.
\end{lemma}

\begin{proof}
    $(\implies)$ Assume $0 \notin \Int P$. Then, we have two cases:
    \begin{enumerate}
        \item If $0 \notin P$, \Cref{thm:sep-hyperplane} gives us a separating hyperplane that contains $F$ in one of the resulting halfspaces. If necessary, take the parallel hyperplane that contains 0.
        \item $0 \in \bd P$, so by \Cref{cor:boundary-face}, $0$ lies in some face of $P$, say $P \cap H$. In which case, $H$ is the desired hyperplane.
    \end{enumerate}

    $(\impliedby)$ Without loss of generality, assume there is a hyperplane, $H$, passing through $0$ such that $F \subseteq H^+$. Since halfspaces are convex, $\conv F \subseteq H^+$. Thus $\conv F$ cannot contain a neighborhood of the origin, that is, $0 \notin \Int P$.
\end{proof}

We now prove that we can characterize conic frames by inspecting the interior of their convex hull.

\begin{theorem}\label{thm:conic-frame-equivalence-farkas}
    Let $F$ be a frame in $\RR^n$ and let $P := \conv F$. $0 \in \Int P$ if and only if we can write any vector in $\R^n$ as a conic combination of the vectors that define $F$.
\end{theorem}

In the proof we will use Farkas' Lemma, a well-known result in optimization:

\begin{lemma}[Farkas' Lemma \citep{ziegler_lectures_1995}]\label{lemma:farkas}
    Let $A \in \R^{n \times m}$ and $z \in \R^n$. Then exactly one of the following is true:
    \begin{enumerate}
        \item There exists a point $x \in \R^m$ with $Ax = z$ and $x \ge 0$
        \item There exists a vector $d \in \R^n$ with $A^\top d \ge 0$ and $d^\top z < 0$
    \end{enumerate}
\end{lemma}

\begin{proof}[Proof of \Cref{thm:conic-frame-equivalence-farkas}]
    We proceed by proving the negation of each statement implies the other.
    
    $(\implies)$ Assume that there is a vector $z \in \R^n$ that cannot be written as a conic combination of the vectors of $F$, so assertion 1 of \Cref{lemma:farkas} is false. Then assertion 2 of \Cref{lemma:farkas} must be true, and thus there is a $d \in \R^n$ with $F^\top d \ge 0$. Let $H := \{ x \in \R^n \ | \ \ll x, d \rr = 0 \}$. Then $0 \in H$ and $F \subseteq H^+$, so by \Cref{lemma:origin-hyperplane}, $0 \notin \Int P$.

    $(\impliedby)$ Now assume $0 \notin \Int P$, so by \Cref{lemma:origin-hyperplane} there is a hyperplane passing through the origin such that $F$ is contained in one of its halfspaces. That is, there is a unit vector $u \in \R^n$ with $\ll f_i, u \rr \ge 0$ for all $i$ and thus $F^\top u \ge 0$. But then assertion 2 of \Cref{lemma:farkas} is true for $z = -u$ since $(-u)^\top u = -1 < 0$, so assertion 1 of \Cref{lemma:farkas} must be false, i.e. $-u$ cannot be written as a conic combination of the vectors that define $F$.
\end{proof}

This proof is fairly concise, but it does not provide much intuition for why it is true. We thus present a different proof for the forward direction which allows us to say something stronger:

\begin{theorem}
\label{thm:conic-frame-equivalence-geometric}
Let $F=\{f_1,\dots,f_m\}$ be a frame in $\R^n$,
and let $P=\conv F$. Then $0\in\Int P$ if and
only if every vector in $\R^n$ is a conic
combination of at most $n$ vectors of $F$.
\end{theorem}

The proof of this result relies on the fact that we can scale $\conv F$ so that any point lies on its boundary, and thus one of the faces of $\conv F$. This provides a scaling of a convex combination which results in a conic combination. We formalize this scaling in the following lemma:

\begin{lemma}[Scaling Lemma]\label{lemma:scaling}
    Let $K$ be a compact convex set in $\RR^n$ such that $0 \in \Int K$ and let $x \neq 0$ be a point in $\RR^n$. Then there exists an $r > 0$ such that $x$ lies on the boundary of $rK$.
\end{lemma}

\begin{proof}
Let
\(
T:=\{t\ge0:tx\in K\}.
\)
Since \(K\) is closed, \(T\) is closed. Boundedness of \(K\) gives an \(R>0\) such that
\(
t\in T\implies t\|x\|\le R,
\)
so \(T\) is bounded. Moreover, \(0\in\operatorname{int}K\) implies that \(T\) contains a positive number. Hence \(M:=\max T\) exists and is positive.

By maximality, \(Mx\in K\), while \((M+h)x\notin K\) for every \(h>0\). These outside points approach \(Mx\), so \(Mx\in\partial K\). Taking \(r=1/M\), positive scaling preserves boundaries and therefore
\(
x=r(Mx)\in r\,\partial K=\partial(rK).
\)
\end{proof}

We now prove \Cref{thm:conic-frame-equivalence-geometric}:

\begin{proof}
$(\implies)$ Suppose $0\in\Int P$, and let
$x\in\R^n$. If $x=0$, use $0=0f_1$.
Now assume $x\neq0$.

By \Cref{lemma:scaling}, there exists $r>0$
such that $x\in\partial(rP)$. By
\Cref{cor:boundary-face,cor:faces-n-1},
$x$ belongs to a proper face $G$ of $rP$
of dimension $k\le n-1$, where
\(
G=\conv\{rf_i:i\in I\}
\)
for some $I\subseteq\{1,\dots,m\}$.

Applying Carath\'eodory's theorem in the
$k$-dimensional affine hull of $G$, there exist
indices $i_1,\dots,i_q\in I$, with
$q\le k+1\le n$, and coefficients
$\lambda_1,\dots,\lambda_q\ge0$ such that
\[
x=\sum_{j=1}^{q}\lambda_j(rf_{i_j}),
\qquad
\sum_{j=1}^{q}\lambda_j=1.
\]
Therefore,
\[
x=\sum_{j=1}^{q}(r\lambda_j)f_{i_j}
\]
is a conic combination of at most $n$ vectors
of $F$.

$(\impliedby)$ Suppose every vector in $\R^n$
is a conic combination of at most $n$ vectors
of $F$. Then $F$ is a conic frame, so
$0\in\Int P$ by
\Cref{thm:conic-frame-equivalence-farkas}.
\end{proof}

Taking \Cref{lemma:origin-strict-pos} and \Cref{thm:conic-frame-equivalence-geometric} together, we have the following:

\begin{corollary}\label{cor:strict-pos-iff-conic}
    Let $F$ be a frame in $\R^n$ consisting of $m$ vectors and let $P = \conv F$. We can write any vector in $\R^n$ as the conic combination of at most $n$ of the vectors that define $F$ if and only if 0 can be written as a strictly positive convex combination of the vectors of $F$.
\end{corollary}

\Cref{cor:strict-pos-iff-conic} gives us a simple method to check if a frame is conic: we just need to check if $0$ can be represented by a strictly positive convex combination of the frame vectors. This is easily verified by solving the following feasibility problem:

\begin{align*}
    \text{find} \quad & x \\
    \text{s.t.} \quad & Fx = 0 \\
    & x > 0 \\
    & \mathbbm{1}^\top x = 1
\end{align*}

We note here that Theorem 3.6 in~\cite{davis1954theory} shows that a frame is conic if and only if 0 can be written as a strictly positive combination of the vectors of the frame. Our proof provides a new perspective and intuition by looking at the polytope defined by a frame.

We summarize the results of this section in the following theorem.
\begin{theorem}\label{thm:conic-frame-equivs}
    Let $F = \{x_1, \dots, x_m\}$ be a frame in $\R^n$ and let $P = \conv(F)$. The following are equivalent:
    \begin{enumerate}[label=(\alph*)]
        \item $0 \in \Int P$
        \item $0$ can be written as a strictly positive convex combination of $x_1, \dots, x_m$
        \item There is no hyperplane passing through $0$ such that $F$ is contained in one of its resulting halfspaces
        \item Any vector in $\R^n$ can be written as a conic combination of the vectors that define $F$
        \item Any vector in $\R^n$ can be written as the conic combination of at most $n$ of the vectors that define $F$
    \end{enumerate}
\end{theorem}

\subsection{Conic Frames and Polyhedra}
Recall that a polyhedron, $C$, is the intersection of finitely many closed halfspaces, 
\[
C=\{x\in\mathbb R^n:\langle s_j,x\rangle\le r_j,\ j=1,\ldots,m\},
\]
where the \(s_j\) are nonzero. A bounded polyhedron is in fact a polytope. And as the following theorem shows, there is a one-to-one correspondence between bounded polyhedra and conic frames.

\begin{theorem}
A nonempty polyhedron \(C\) is bounded if and only if its defining normals \(s_1,\ldots,s_m\) form a conic frame for \(\mathbb R^n\).    
\end{theorem}

Before proving this theorem, we introduce the notion of a normal cone as well as a characterization of the normal cone of a polyhedron which can be found in  (\cite{HiriartUrruty_ConvexAnalysis} Example 5.2.6(b)).

\begin{definition}[Normal Cone]
    The direction $s \in \R^n$ is said to be \emph{normal} to $C$ at $x \in C$ when $\ll s, y-x \rr \le 0$ for all $y \in C$. The set of all such directions is called the \emph{normal cone} to $C$ at $x$, denoted by $N_C(x)$.
\end{definition}

Using the definition of a polyhedron above, define $J(x) := \{j = 1,\dots,m : \ll s_j, x \rr = r_j\}$, the index-set of active constraints at $x \in C$. Then $N_C(x) = \{\sum\limits_{j \in J(x)} \lambda_j s_j : \lambda_j \ge 0\}$.

\begin{proof}
($\implies$) Suppose \(C\) is bounded. Because it is an intersection of closed halfspaces, it is closed; hence it is compact.

Take any \(v\in\mathbb R^n\). The continuous linear functional
\(
x\longmapsto\langle v,x\rangle
\)
attains its maximum on \(C\), say at \(x_v\in C\). Therefore, for every \(y\in C\),
\(
\langle v,y\rangle\le\langle v,x_v\rangle,
\)
or equivalently,
\(
\langle v,y-x_v\rangle\le0.
\)
By the definition of the normal cone, this says
\(
v\in N_C(x_v) =\operatorname{cone}\{s_j:j\in J(x_v)\},
\)
where
\(
J(x_v)=\{j:\langle s_j,x_v\rangle=r_j\}.
\)
Consequently, there are nonnegative coefficients such that
\[
v=\sum_{j\in J(x_v)}\lambda_j s_j.
\]
Since \(v\) was arbitrary,
\(
\operatorname{cone}\{s_1,\ldots,s_m\}=\mathbb R^n.
\)
Thus the normals positively span \(\mathbb R^n\). They also span it linearly, so they form a finite frame, and therefore a conic frame.

($\impliedby$) Suppose the normals form a conic frame. For each coordinate vector \(e_i\), choose nonnegative coefficients \(a_{ij},b_{ij}\) such that
\[
e_i=\sum_{j=1}^m a_{ij}s_j,
\qquad
-e_i=\sum_{j=1}^m b_{ij}s_j.
\]
For any \(x\in C\), the defining inequalities then give
\[
x_i=\sum_{j=1}^m a_{ij}\langle s_j,x\rangle
\le\sum_{j=1}^m a_{ij}r_j,
\]
and similarly,
\[
-x_i\le\sum_{j=1}^m b_{ij}r_j.
\]
Hence
\[
-\sum_{j=1}^m b_{ij}r_j
\le x_i\le
\sum_{j=1}^m a_{ij}r_j
\qquad(i=1,\ldots,n).
\]
Every coordinate is bounded above and below by fixed finite constants, so \(C\) is bounded.

\end{proof}

\subsection{Random Conic Frames}
\label{subsec:random-method}

We now turn to randomly generating a conic frame with high probability. To do so, we begin by stating Wendel's Theorem, which provides the probability that a set of randomly generated vectors on the unit sphere all lie on some hemisphere.

\begin{theorem}[Wendel's Theorem~\cite{wendel1962problem}]\label{thm:wendel}
    Let $x_1, x_2, \dots x_m$ be random unit vectors in $\R^n$ whose joint distribution is invariant under all reflections through the origin, i.e., independently changing the signs of any subset of the vectors, and is such that with probability one all subsets of size $n$ are linearly independent, then the probability, $p_{n,m}$, that all the vectors lie on some hemisphere is:
    $$ p_{n,m} = 2^{-m+1} \sum\limits_{k=0}^{n-1} \binom{m-1}{k}$$
\end{theorem}

Note that $\binom{m-1}{k} = 0$ for $m-1 < k$, so $p_{n,m} = 1$ for $m \le n$.

A hemisphere is determined by a hyperplane through the origin, so $p_{n, m}$ is equivalently the probability that there exists a hyperplane through the origin defining a halfspace that contains $x_1, \dots, x_m$.

\begin{theorem}[Probability that a random frame is conic]\label{thm:prob-conic}
    Let $F$ be a frame consisting of $m$ vectors in in $\R^n$ satisfying the hypotheses of \Cref{thm:wendel}, then $F$ is a conic frame with probability $1 - p_{n,m}$.
\end{theorem}

\begin{proof}
    By \Cref{thm:conic-frame-equivs}, $F$ is conic if and only if there is no hyperplane through the origin defining a halfspace that contains $F$, which by \Cref{thm:wendel} is precisely $1 - p_{n,m}$.
\end{proof}

We note that a set of Gaussian vectors satisfies the hypotheses of \Cref{thm:prob-conic} almost surely.

\bigskip
A natural question that arises is whether we can turn a given frame easily into a conic frame? There are at least  two trivial ways to turn any frame $\{f_k\}_{k=1}^m$ into a conic frame: 1. Just add all negative elements $\{-f_k\}_{k=1}^m$. While this preserves the condition number of the frame operator, it doubles the number of elements in the frame, which may not be desirable. 2. Add $f_{m+1}:= -\frac{1}{m}\sum_{k=1}^m f_k$ to the frame. A simple calculation shows that the resulting frame $\{f_k\}_{k=1}^{m+1}$ is conic.
While this approach adds only one element, the condition number of the frame operator of this so augmented frame is worse.

A more interesting way to turn a frame into a conic frame is by switching the signs of some frame elements. This procedure neither increases the number of elements, nor does it change the frame bounds. But can switching some signs of a frame always yield a conic frame? Moreover,
without further assumptions about the frame, a deterministic search through all signs is computationally not feasible if $m$ is large.
Is there an efficient randomized algorithm to do so?

For specific cases, we will provide both a deterministic and a randomized algorithm to turn a frame into a conic frame via sign changes.

The following theorem  provides a method to turn a frame into a conic frame with high probability.
We need only ensure that we apply a distribution to the frame vectors that satisfies the hypotheses for \Cref{thm:prob-conic}.

\begin{theorem}\label{thm:make-conic}
    Let $F$ be a frame in $\R^n$ such that every subset of $n$ of its columns is linearly independent and let $B_1, B_2, \dots, B_m$ be i.i.d. random variables such that $\P(B_i = 1) = \frac{1}{2}$ and $\P(B_i = -1) = \frac{1}{2}$. If $b_i$ is a realization of $B_i$ and $B$ is the $m \times m$ diagonal matrix with $B_{ii} = b_i$, and $G$ is the matrix with columns equal to the unit normalized columns of $F$ then $GB$ is a conic frame with probability $1 - p_{n,m}$.
\end{theorem}

\begin{proof}
    First note that normalization preserves conicity and linear independence. Then note that the joint distribution of the columns of $GB$ are invariant under all reflections through the origin by construction. Thus $FB$ satisfies the conditions of \Cref{thm:wendel}, so the result follows by \Cref{thm:prob-conic}.
\end{proof}

The reader may wonder how large does $m$ have to be compared to $n$ for this method to have a reasonable probability of success? The probability of success is exactly 1/2 when $m=2n$. 
It is not difficult to see that by Hoeffding's inequality, $1-p_{n,m}\ge1-\delta$ whenever
\[
m\ge2n+2\sqrt{n\log(1/\delta)}+2\log(1/\delta),
\]
while for $m \le 2n-1$ the probability of success becomes exponentially small.

Regarding dictionary learning, given a full spark~\cite{donoho_optimally_2003,alexeev2012fullspark}) $n$-dimensional dictionary, i.e. every size $n$ subset of its atoms is linearly independent, \Cref{thm:make-conic} can be used to make it conic. This method leaves the frame bound and sparsity of the system unchanged.

\subsection{A Deterministic Construction by Sign Changes}
\label{subsec:deterministic-method}
In addition to our randomized method, we also present a deterministic method to make a full spark frame conic.

\begin{theorem}
     Let \(F=[f_1,\ldots,f_m]\in\mathbb R^{n\times m}\) be full spark, with \(m\ge n+1\). A diagonal matrix \(D\) with diagonal entries in \(\{-1,1\}\) can be constructed deterministically such that \(FD\) is a conic frame.
\end{theorem}

\begin{proof}
Since \(F\) is full spark, \(A=[f_1,\ldots,f_n]\) is invertible. Solve
\[
Ac=f_{n+1}.
\]
Every \(c_i\) is nonzero: otherwise \(f_{n+1}\) would lie in the span of \(n-1\) of the first \(n\) columns, contradicting full spark.

Define
\[
g_i=\operatorname{sgn}(c_i)f_i\quad(1\le i\le n),
\qquad
g_{n+1}=-f_{n+1}.
\]
Then
\[
\sum_{i=1}^{n}|c_i|g_i+g_{n+1}=0.
\]
Thus these \(n+1\) vectors span \(\mathbb R^n\) and admit a dependence with strictly positive coefficients. so by \Cref{thm:conic-frame-equivs}, they positively span \(\mathbb R^n\). Assign arbitrary signs to the remaining columns; adding vectors preserves positive spanning.    
\end{proof}

This construction also preserves the frame bounds and sparsity and has complexity \(O(n^3)\).

\section{Maximization of Cosine Measure}
\label{sec:max-cosine-measure}

Amidst studying conic frames and the cosine measure, we noticed the similarity between coherence for a frame and cosine measure. Noting that ETFs achieve minimal coherence, it seems obvious to conjecture that ETFs achieve maximal cosine measure. However, one major difference between cosine measure and coherence is that the formula for cosine measure does not include an absolute value. A natural fix for this is to consider instead a Twin-ETF, defined below.

\begin{definition}
    Let $F$ be an Equiangular Tight Frame. Then the set of vectors in $F$ and $-F$, denoted $[F, -F]$, is a Twin-ETF.
\end{definition}

Note that for a Twin ETF, \[
\max_{g\in[F,-F]}\langle v,g\rangle
=\max_{f\in F}|\langle v,f\rangle|
\]
for every $v$.

The other difference between minimal coherence and cosine measure is that cosine measure checks each vector in the underlying space, while coherence only looks at the frame vectors themselves.
This lead us to conjecture that perhaps regular polytopes (polytopes with maximal symmetry) might maximize cosine measure. In particular, it seems intuitively clear that a regular polytope that is also a twin-ETF, an object that is optimal for several geometrically inspired measures, should also be optimal with respect to the cosine measure.
Alas, we were surprised to find that neither of these conjectures turn out to be true. In fact, a counterexample for 8 vectors in $\R^3$ demonstrates this.

Begin by noting that the following frame is a unit norm ETF in $\R^3$.

$$F = \frac{1}{\sqrt{3}}\begin{bmatrix}
    1 & 1 & 1 & 1\\
    1 & 1 & -1 & -1\\
    1 & -1 & 1 & -1\\
\end{bmatrix}$$

Thus $[F, -F]$ is a Twin ETF. Additionally, it is not difficult to see that $[F, -F]$ is in fact the unit norm cube and is thus also a regular polytope. However, we will demonstrate shortly that this conic frame does not have the maximal cosine measure for 8 vector frames in $\R^3$. In hindsight this is less surprising, since the coherence of a frame is a {\em packing} quantity, while the cosine measure is a {\em covering} quantity.

The square-antiprism is similar to a cube, but the top and bottom faces are rotated $45^\circ$ from each other. That is, for some value of the parameter $0 < a < 1$ with $r = \sqrt{1 - a^2}$, a square antiprism has the following form:

$$S = \begin{bmatrix}
r\cos(0) & r\cos(\frac{\pi}{2}) & r\cos(\pi) & r\cos(\frac{3\pi}{2}) & r\cos(\frac{\pi}{4}) & r\cos(\frac{3\pi}{4}) & r\cos(\frac{5\pi}{4}) & r\cos(\frac{7\pi}{4}) \\
r\sin(0) & r\sin(\frac{\pi}{2}) & r\sin(\pi) & r\sin(\frac{3\pi}{2}) & r\sin(\frac{\pi}{4}) & r\sin(\frac{3\pi}{4}) & r\sin(\frac{5\pi}{4}) & r\sin(\frac{7\pi}{4}) \\
a        & a                    & a          & a                     & -a                   & -a                    & -a                    & -a
\end{bmatrix}$$

For $a = 0.7$, we have the following:

$$S = \begin{bmatrix}
    0.7141 & 0      & -0.7141 & 0       &  0.505 & -0.505 & -0.505 &  0.505 \\
     0 & 0.7141 &  0      & -0.7141 &  0.505 &  0.505 & -0.505 & -0.505 \\
0.7    & 0.7    &  0.7    &  0.7    & -0.7   & -0.7   & -0.7   & -0.7
\end{bmatrix}$$

We implemented and ran the deterministic method from~\cite{hare2020} in MATLAB to calculate the cosine measure of each of these conic frames. We obtained $\cm([F, -F]) \approx 0.5774$ while $\cm(S) \approx 0.6029$, meaning that this square antiprism has higher cosine measure than the unit norm cube. At first glance, this is surprising. However, in~\cite{dodangeh_optimal_2016} it is shown that finding the set of $m$ vectors in $\R^n$ with maximal cosine measure is equivalent to finding the most economical covering by $S^{n-1}$ by $m$ spherical caps. This problem is similar to (but still different from) the Tammes problem which looks for the spherical code of $m$ points with largest angular separation between points. It is known that the solution for $8$ points in $\R^3$ is the square antiprism~\cite{conway_sphere_1999}.

But this is not the whole story. The best known 8-cap covering of the sphere is a $D_2d$ point configuration~\cite{HardinSloaneSmith1994coverings} of the following form\footnote{As an abstract group, $D_2d$ is isomorphic to the dihedral group of order 8.}:
\begin{gather*} F=\bigl\{(\pm s_1,0,c_1),\ (0,\pm s_1,-c_1),\ (0,\pm s_2,c_2),\ (\pm s_2,0,-c_2)\bigr\},\\ c_1\approx0.7867473615,\quad c_2\approx 0.1715313650,\quad s_i=\sqrt{1-c_i^2}, \end{gather*}
which has $\operatorname{cm}(F)\approx0.6673188866$, an even higher cosine measure. Also note that the Tammes-optimal square anti-prism~\cite{dostert_exact_sdp_2020} has cosine measure $a \approx 0.511$, demonstrating that these problems are distinct.

\begin{section}{The Cosine Measure of Gaussian Frames}
\label{sec:cm-gaussianframes}

\subsection{Facet Height and Cosine Measure}

We now provide an alternative characterization of the cosine measure of a polytope, namely the minimum facet height of the polytope.

\begin{lemma}[Cosine measure as the minimum facet height]
\label{lem:cosine-minimum-facet}
Let
\[
X=\{x_1,\dots,x_m\}\subset\mathbb{S}^{d-1}
\]
and let
\[
P:=\operatorname{conv}(X).
\]
Suppose that
\[
0\in\operatorname{int}P.
\]
For each facet $F\in\mathcal{F}(P)$, let $u_F\in\mathbb{S}^{d-1}$
be its outward unit normal and let $H_F>0$ be its height, so that
\[
F\subseteq\{x\in\mathbb{R}^d:\langle u_F,x\rangle=H_F\}
\]
and
\[
P\subseteq\{x\in\mathbb{R}^d:\langle u_F,x\rangle\le H_F\}.
\]
Then
\[
\cm X
=
\min_{\|u\|=1}\max_{1\le i\le m}\langle u,x_i\rangle
=
\min_{F\in\mathcal{F}(P)}H_F.
\]
Equivalently, the cosine measure is the radius of the largest
Euclidean ball centered at the origin and contained in $P$.
\end{lemma}

\begin{proof}
Let
\[
h_P(u):=\max_{x\in P}\langle u,x\rangle,
\qquad u\in\mathbb{S}^{d-1},
\]
be the support function of $P$. Because $P$ is the convex hull of
$x_1,\dots,x_m$, a linear functional attains its maximum over $P$ at
one of these points. Therefore,
\[
h_P(u)=\max_{1\le i\le m}\langle u,x_i\rangle.
\]
It follows that
\begin{equation}
\cm X
=
\min_{\|u\|=1}h_P(u).
\label{eqn:cm-support}
\end{equation}

Let
\[
r_0(P):=
\sup\left\{
r\ge0:r\mathbb{B}_2^d\subseteq P
\right\}
\]
denote the radius of the largest Euclidean ball centered at the
origin and contained in $P$. We first show that
\begin{equation}
r_0(P)=\min_{\|u\|=1}h_P(u).
\label{eqn:large-ball-support}
\end{equation}

For every $r\ge0$, we claim that
\begin{equation}
r\mathbb{B}_2^d\subseteq P
\quad\Longleftrightarrow\quad
r\le h_P(u)
\quad\text{for every }u\in\mathbb{S}^{d-1}.
\label{eqn:ball-support}
\end{equation}

Suppose first that $r\mathbb{B}_2^d\subseteq P$. The support function
of $r\mathbb{B}_2^d$ is
\[
h_{r\mathbb{B}_2^d}(u)
=
\max_{\|x\|\le r}\langle u,x\rangle
=
r,
\qquad \|u\|=1.
\]
By monotonicity of support functions under set inclusion,
\[
r=h_{r\mathbb{B}_2^d}(u)\le h_P(u)
\]
for every $u\in\mathbb{S}^{d-1}$.

Conversely, suppose that
\[
r\le h_P(u)
\qquad\text{for every }u\in\mathbb{S}^{d-1},
\]
but that $r\mathbb{B}_2^d\not\subseteq P$. Then there exists some
$x\in r\mathbb{B}_2^d$ such that $x\notin P$. Since $P$ is closed and
convex, the strict separation theorem gives a unit vector
$u\in\mathbb{S}^{d-1}$ satisfying
\[
\langle u,x\rangle
>
\sup_{y\in P}\langle u,y\rangle
=
h_P(u).
\]
On the other hand, because $\|x\|\le r$, the Cauchy--Schwarz
inequality gives
\[
\langle u,x\rangle
\le \|u\|\,\|x\|
\le r
\le h_P(u),
\]
which is a contradiction. This proves \eqref{eqn:ball-support}. Taking the largest
admissible value of $r$ in \eqref{eqn:ball-support} gives \eqref{eqn:large-ball-support}.

It remains to express the centered inradius in terms of the facet
heights. The facet half-space representation of $P$ is
\begin{equation}
P
=
\bigcap_{F\in\mathcal{F}(P)}
\left\{
x\in\mathbb{R}^d:
\langle u_F,x\rangle\le H_F
\right\}.
\label{eqn:facet-rep}
\end{equation}
For every facet $F$,
\[
\max_{\|x\|\le r}\langle u_F,x\rangle=r,
\]
because $\|u_F\|=1$. Consequently, \eqref{eqn:facet-rep} implies
\[
\begin{aligned}
r\mathbb{B}_2^d\subseteq P
&\quad\Longleftrightarrow\quad
r\mathbb{B}_2^d
\subseteq
\left\{
x:\langle u_F,x\rangle\le H_F
\right\}
\quad\text{for every }F\\
&\quad\Longleftrightarrow\quad
r\le H_F
\quad\text{for every }F\in\mathcal{F}(P).
\end{aligned}
\]
Therefore,
\begin{equation}
r_0(P)=\min_{F\in\mathcal{F}(P)}H_F.
\label{eqn:ball-facet}
\end{equation}

Combining \eqref{eqn:cm-support}, \eqref{eqn:large-ball-support}, and \eqref{eqn:ball-facet}, we conclude that
\[
\cm X
=
\min_{\|u\|=1}h_P(u)
=
r_0(P)
=
\min_{F\in\mathcal{F}(P)}H_F.
\]
\end{proof}

\subsection{Concentration of the Cosine Measure of Normalized Gaussian Frames}

Besides providing interesting geometric intuition, the facet height characterization connects cosine measure to previous work by Bonnet and O'Reilly~\cite{BonnetOReilly2022} on asymptotic properties of convex hulls of normalized Gaussian vectors. We adapt the
facet-count estimates of
\cite[Theorem~3.4 and its proof in Section~6.2.1]{BonnetOReilly2022}
to obtain quantitative bounds for the cosine measure.
The lower-tail argument follows their Lemma~6.10, with a sharper
binomial-coefficient bound.

Our contributions are the following:
\begin{enumerate}
    \item We demonstrate semi-asymptotic exponential decay of the tail probabilities for sufficiently large dimension.
    \item We provide a new chain of inequalities to correct an inaccuracy in their proof. Note: The conclusion of~\cite[Theorem~3.4]{BonnetOReilly2022} remains valid in the regimes of Theorem 5.2; only the stated bound in~\cite[Lemma~6.10]{BonnetOReilly2022} fails.
\end{enumerate}

\begin{theorem}[Exponential concentration of the cosine measure]
\label{thm:gaussian-cosine-concentration}
Let $g_1,\dots,g_m$ be independent standard Gaussian vectors in
$\mathbb{R}^n$, and define
\[
u_i:=\frac{g_i}{\|g_i\|},
\qquad
P_{n,m}:=\operatorname{conv}\{u_1,\dots,u_m\}.
\]
The vectors $u_1,\dots,u_m$ are independent and uniformly distributed on
$\mathbb{S}^{n-1}$.
Suppose that
\[
n\to\infty
\qquad\text{and}\qquad
\frac{m}{n}\to\infty.
\]
Then the following assertions hold for their cosine measure
\[
\mu_{n,m}
=
\min_{\|v\|=1}\max_{1\le i\le m}\langle u_i, v\rangle.
\]

(Note we use $\mu_{n,m}$ to denote the cosine measure in order to emphasize the dependence on $m$ and $n$.)

\begin{enumerate}
\item \textbf{Subexponential regime.} If $\log m=o(n)$, then, for every fixed $\varepsilon>0$, there
exist constants $c,C>0$ such that
\[
\mathbb{P}\left(
\left|
\frac{\mu_{n,m}}
{\sqrt{2\log(m/n)/n}}
-1
\right|>\varepsilon
\right)
\le Ce^{-cn}
\]
for all sufficiently large $n$.

\item \textbf{Exponential regime.} If
\[
\frac{\log m}{n}\longrightarrow\gamma\in(0,\infty),
\]
then, for every fixed $\varepsilon>0$, there exist constants $c,C>0$
such that
\[
\mathbb{P}\left(
\left|
\mu_{n,m}-\sqrt{1-e^{-2\gamma}}
\right|>\varepsilon
\right)
\le Ce^{-cn}
\]
for all sufficiently large $n$.

\item \textbf{Superexponential regime.} If
\[
\frac{\log m}{n}\longrightarrow\infty,
\]
then, for every fixed $\varepsilon\in(0,1)$, there exist constants
$c,C>0$ such that
\[
\mathbb{P}\left(\mu_{n,m}<1-\varepsilon\right)
\le C\exp\bigl(-c\,n\log m\bigr)
\]
for all sufficiently large $m,n$.
\end{enumerate}

\end{theorem}

\begin{proof}
For $n\ge2$ and $m\ge n+1$,
$P_{n,m}$ is full-dimensional almost surely.
Following the signed-height convention of
\cite[Section~1]{BonnetOReilly2022}, for every facet \(F\),
define its signed height \(H_F\in[-1,1]\) by its outward
unit normal \(u_F\):
\[
F\subseteq\{x:\langle u_F,x\rangle=H_F\},
\qquad
P_{n,m}\subseteq\{x:\langle u_F,x\rangle\le H_F\}.
\]
This definition applies whether or not the event
\[
E_{n,m}:=
\left\{0\in\operatorname{int}P_{n,m}\right\}
\]
occurs. 

\medskip
\noindent
By Lemma~\ref{lem:cosine-minimum-facet}, on the event $E_{n,m}$,
\begin{equation}
\mu_{n,m}
=
\min_{F\in\mathcal{F}(P_{n,m})}H_F,
\label{eqn:cm-facet-thm}
\end{equation}
where $H_F>0$ on $E_{n,m}$.

\medskip
\noindent

Define
\[
h_1
=
\sqrt{
1-
\left(
\frac{r_1n(\log(m/n))^{3/2}}{m}
\right)^{\frac{2}{n-1}}
}
\]
and
\begin{equation}
h_2
=
\sqrt{
1-
\left(
\frac{r_2n}{m}
\right)^{\frac{2(n+1)}{(n-1)^2}}
},
\end{equation}
where $r_1>0$ is sufficiently large and $r_2>0$ is sufficiently small;
see \cite[Theorem~3.4 and equation~(6.24)]{BonnetOReilly2022}.
Their notation uses $n$ in place of our $m$, and $d$ in place of our
$n$. For all sufficiently large $n$, $0<h_1<h_2<1$.

Let
\[
N_-:=
\#\{F:H_F\in[-1,h_1]\}
\]
and
\[
N_+:=
\#\{F:H_F\in[h_2,1]\}.
\]

By the definition of $F_{[a,b]}$ in \cite{BonnetOReilly2022},
\begin{equation}
\mathbb{E}N_-=F_{[-1,h_1]},
\qquad
\mathbb{E}N_+=F_{[h_2,1]}.
\end{equation}

We follow the facet-count argument in the proof of
\cite[Theorem~3.4]{BonnetOReilly2022}, particularly the proofs
of Lemmas~6.10 and~6.11. The strategy is to bound the expected
number of facets whose heights fall outside the interval
\([h_1,h_2]\). Markov's inequality then controls the probability
that any such facet exists. On the event that the origin lies
in the interior of the polytope, we use~\Cref{lem:cosine-minimum-facet}
to transfer this simultaneous control of facet heights to control of the cosine measure.

Repeating the proof of \cite[Lemma~6.10]{BonnetOReilly2022}
with a sharper binomial-coefficient estimate gives the following
bound. For $\beta>-1$, let
\[
c_\beta
:=
\left(\int_{-1}^{1}(1-t^2)^\beta\,dt\right)^{-1}.
\]
Thus
\[
c_\beta\int_{-1}^{1}(1-t^2)^\beta\,dt=1.
\]

And set
\[
L:=\log(m/n),
\qquad
A:=c_{\frac{n-3}{2}}
\int_{h_1}^{1}(1-t^2)^{\frac{n-3}{2}}\,dt.
\]
The estimate from
\cite[Lemma~6.9]{BonnetOReilly2022} with \(D=(n-3)/2\),
together with the definition of $h_1$, gives,
\[
A\ge
\frac{c_{(n-3)/2}}{n-1}
\frac{(1-h_1^2)^{(n-1)/2}}{h_1}
\left(1-\frac{1-h_1^2}{(n+1)h_1^2}\right).
\]
By \cite[equation~(6.2)]{BonnetOReilly2022},
\[
\frac{c_{(n-3)/2}}{n-1}
\sim \frac{1}{\sqrt{2\pi n}}.
\]

Moreover, \cite[Lemma~6.8]{BonnetOReilly2022} gives
\(nh_1^2\to\infty\), so the factor in parentheses tends to \(1\). Absorbing fixed numerical factors into \(b_0>0\) gives:

\[
A\ge b_0\frac{(1-h_1^2)^{(n-1)/2}}{h_1\sqrt n}
\]

Following the estimate in
\cite[equation~(6.25)]{BonnetOReilly2022}, write
\[
h_1^2=1-e^{-x},
\qquad
x=\frac{2}{n-1}
\log\left(\frac{m}{r_1nL^{3/2}}\right).
\]
Eventually \(x>0\) and \(r_1L^{3/2}\ge1\). Using \(1-e^{-x}\le x\),
\[
h_1^2
\le \frac{2}{n-1}
\log\left(\frac{m}{r_1nL^{3/2}}\right)
\le\frac{2L}{n-1}.
\]
Therefore, for \(n\ge2\),
\[
h_1\sqrt n
\le\sqrt{\frac{2nL}{n-1}}
\le2\sqrt L.
\]
Substituting the definition of $h_1$ and this denominator bound gives
\[
b_0\frac{(1-h_1^2)^{(n-1)/2}}{h_1\sqrt n}
\ge
b_0\frac{r_1nL^{3/2}/m}{2\sqrt L}
=\frac{b_0}{2}\frac{r_1nL}{m}.
\]
Setting $b_1 = b_0 / 2$, we have,
\[
A\ge b_0\frac{(1-h_1^2)^{(n-1)/2}}{h_1\sqrt n}
\ge b_1\frac{r_1nL}{m}.
\]
Using our notation, \cite[equation~(6.4)]{BonnetOReilly2022} states

\[
\begin{aligned}
F_{[a,b]}
={}&2\binom{m}{n}\,
c_{\frac{n^2-2n-1}{2}}
\int_a^b
(1-h^2)^{\frac{n^2-2n-1}{2}}
\left(
c_{\frac{n-3}{2}}
\int_{-1}^{h}(1-s^2)^{\frac{n-3}{2}}\,ds
\right)^{m-n}\,dh.
\end{aligned}
\]
The following chain of inequalities corrects the inaccuracy explained in \Cref{sec:explain-inacc}. It is inspired by~\cite[Lemma~6.11]{BonnetOReilly2022}. Bounding the inner normalized integral by $1-A$
and the remaining normalized integral by $1$
gives the first inequality below.
Using $\binom mn\le(em/n)^n$, \(1-A\le e^{-A}\), and $m-n\ge m/2$ eventually,
\[
\begin{aligned}
\mathbb{E}N_-
&\le 2\binom mn(1-A)^{m-n}\\
&\le 2\left(\frac{em}{n}\right)^n
   \exp\bigl(-b r_1nL\bigr)\\
&=2e^n\left(\frac mn\right)^{-(br_1-1)n},
\end{aligned}
\]
where $b = b_1 / 2>0$ is an absolute constant. Choose \(r_1\) so that \(c:=br_1-1>0\). Because \(L\to\infty\), eventually
\begin{equation}
\mathbb EN_-
\le e^{\log2+n-cnL}
\le e^{-(c/2)nL}
=\left(\frac mn\right)^{-an},
\qquad a:=c/2>0.
\label{eqn:en-minus-bound}
\end{equation}

By \cite[Lemma~6.11, equation~(6.27)]{BonnetOReilly2022},
\begin{equation}
\mathbb{E}N_+
\le
\frac{C\sqrt n}{me}(r_2e)^{n+1}.
\label{eqn:en-plus-bound}
\end{equation}
We may choose $r_2$ sufficiently small that $r_2e<1$. Setting
\[
c_2:=-\log(r_2e)>0,
\]
equation \eqref{eqn:en-plus-bound} gives
\begin{equation}
\mathbb{E}N_+
\le
C\frac{\sqrt n}{m}e^{-c_2n},
\label{eqn:en-plus-bound-final}
\end{equation}
after changing the constant $C$.

\medskip
\noindent
\textbf{Probability bounds.}
On the event $E_{n,m}$, equation \eqref{eqn:cm-facet-thm} implies that $\mu_{n,m}<h_1$
forces $N_->0$, while $\mu_{n,m}>h_2$ forces $N_+>0$.
On $E_{n,m}^{\,c}$, the facet half-space representation gives a facet
with $H_F\le0<h_1$, so $N_->0$ as well. Hence, almost surely,
\begin{equation}
\{\mu_{n,m}\notin[h_1,h_2]\}
\subseteq
\{N_->0\}
\cup\{N_+>0\}.
\label{eqn:set-inclusion}
\end{equation}
Since \(E_{n,m}^{c}\subseteq\{N_->0\}\), the lower
facet-count estimate also controls the probability that the
frame fails to positively span \(\mathbb R^n\).
The union bound and Markov's inequality therefore give
\begin{equation}
\begin{aligned}
\mathbb{P}\bigl(\mu_{n,m}\notin[h_1,h_2]\bigr)
&\le
\mathbb{P}(N_->0)
+\mathbb{P}(N_+>0)\\
&\le
\mathbb{E}N_-
+\mathbb{E}N_+.
\end{aligned}
\label{eqn:markov}
\end{equation}

Combining \eqref{eqn:en-minus-bound}, \eqref{eqn:en-plus-bound-final}, and \eqref{eqn:markov}, we obtain the
quantitative cosine-measure estimate
\begin{equation}
\mathbb{P}\bigl(\mu_{n,m}\notin[h_1,h_2]\bigr)
\le
e^{-an\log(m/n)}
+
C\frac{\sqrt n}{m}e^{-c_2n}.
\label{eqn:cm-combined}
\end{equation}
Since $m/n\to\infty$, eventually $\log(m/n)\ge1$, and hence
\[
e^{-an\log(m/n)}\le e^{-an}.
\]
The right-hand side of \eqref{eqn:cm-combined} is therefore at most
\begin{equation}
C'e^{-c'n}
\label{eqn:cm-final}
\end{equation}
for suitable constants $c',C'>0$.

It remains to identify the deterministic asymptotics of the
cutoffs \(h_1\) and \(h_2\). These follow from
\cite[Lemma~6.8]{BonnetOReilly2022}, applied to the cutoffs
as explained in their proof of Theorem~2.4.

\medskip
\noindent
\textbf{Subexponential regime.}
By \cite[Lemma~6.8(i) and the proof of Theorem~2.4]{BonnetOReilly2022},
\begin{equation}
\frac{h_i}{\sqrt{2\log(m/n)/n}}\longrightarrow1,
\qquad i=1,2.
\end{equation}
Thus, for every fixed $\varepsilon>0$, eventually
\[
[h_1,h_2]
\subseteq
\left[
(1-\varepsilon)\sqrt{\frac{2\log(m/n)}n},
(1+\varepsilon)\sqrt{\frac{2\log(m/n)}n}
\right].
\]
Combining this deterministic inclusion with \eqref{eqn:cm-final} proves part~(1).

\medskip
\noindent
\textbf{Exponential regime.}
By \cite[Lemma~6.8(ii)]{BonnetOReilly2022},
\begin{equation}
h_1,h_2\longrightarrow\sqrt{1-e^{-2\gamma}}.
\end{equation}
Hence, for every fixed $\varepsilon>0$, eventually
\[
[h_1,h_2]
\subseteq
\left[
\sqrt{1-e^{-2\gamma}}-\varepsilon,
\sqrt{1-e^{-2\gamma}}+\varepsilon
\right].
\]
Combining this with \eqref{eqn:cm-final} proves part~(2).

\medskip
\noindent
\textbf{Superexponential regime.}
Lemma~6.8(iii) of \cite{BonnetOReilly2022} gives
\begin{equation}
h_1\longrightarrow1.
\end{equation}
Therefore, for every fixed $\varepsilon\in(0,1)$, eventually
$h_1>1-\varepsilon$. The same argument as in \eqref{eqn:set-inclusion}, followed by
\eqref{eqn:en-minus-bound}, gives
\begin{equation}
\mathbb{P}(\mu_{n,m}<1-\varepsilon)
\le
\mathbb{E}N_-\le e^{-an\log(m/n)}.
\end{equation}
Since $(\log m)/n\to\infty$,
\[
\frac{\log n}{\log m}\to0,
\qquad
\log(m/n)=\log m-\log n\sim\log m.
\]
Consequently, eventually
\[
e^{-an\log(m/n)}\le e^{-(a/2)n\log m}.
\]
This proves part~(3).
\end{proof}

\subsection{Non-asymptotic Gaussian Tail Bounds in the Subexponential Regime}

In the preceding section, we proved semi-asymptotic convergence rates for the cosine measure of Gaussian frames in three regimes.
We now provide non-asymptotic tail bounds for the cosine measure of Gaussian frames in the sub-exponential regime, beginning with a lower bound and then an upper bound. The proof techniques are very different from the ones used in the previous section.

\subsubsection{A lower bound for the Gaussian cosine measure}

Let
$$
g_1,\dots,g_m\overset{\mathrm{iid}}{\sim}N(0,I_n),
\qquad
u_i:=\frac{g_i}{\norm{g_i}},
$$
define
$$
K_{n,m}
:=
\min_{\norm{v}=1}
\max_{1\le i\le m}\ip{g_i}{v},
$$
and note the cosine measure,
$$
\mu_{n,m}
=
\min_{\norm{v}=1}
\max_{1\le i\le m}\ip{u_i}{v}.
$$
Set
$$
\alpha:=\frac mn.
$$

\begin{theorem}
\label{thm:lower}
Fix $\varepsilon\in(0,1)$ and set
$$
\rho:=\frac{\varepsilon}{3},
\qquad
t:=(1-\rho)\sqrt{2\log\alpha},
\qquad
\beta:=1-(1-\rho)^2=2\rho-\rho^2.
$$
There exist absolute constants $C,c>0$ and a number
$\alpha_0(\varepsilon)$ such that, whenever
$\alpha\ge\alpha_0(\varepsilon)$,
\begin{equation}
\Pp\{K_{n,m}<t\}
\le
C\,
\frac{\sqrt{\log\alpha}}
{n\alpha^\beta}
+
2e^{-cn}.
\label{eq:unnormalized-lower-tail}
\end{equation}

Moreover,
\begin{align}
&\Pp\left\{
\mu_{n,m}
<
(1-\varepsilon)
\sqrt{\frac{2\log(m/n)}{n}}
\right\}
\nonumber\\
&\qquad\le
C\,
\frac{\sqrt{\log(m/n)}}
{n(m/n)^\beta}
+
2e^{-cn}
+
m\exp\left(-\frac{n\varepsilon^2}{18}\right),
\label{eq:normalized-lower-tail}
\end{align}
where
$$
\beta
=
\frac{2\varepsilon}{3}
-
\frac{\varepsilon^2}{9}.
$$

In particular, if
$$
\frac mn\to\infty,
\qquad
\log\frac mn=o(n),
$$
then for every fixed $\varepsilon>0$,
$$
\Pp\left\{
\mu_{n,m}
<
(1-\varepsilon)
\sqrt{\frac{2\log(m/n)}{n}}
\right\}
\longrightarrow0.
$$
\end{theorem}

\begin{proof}

\textbf{Step 1: Gaussian min--max formulation.}

Let $G\in\R^{m\times n}$ be the matrix whose $i$th row is $g_i^\top$,
and let
$$
\Delta_m
=
\left\{
y\in\R^m:
y_i\ge0,\quad
\sum_{i=1}^m y_i=1
\right\}.
$$
Since
$$
\max_{1\le i\le m}x_i
=
\max_{y\in\Delta_m}y^\top x,
$$
we have
\begin{equation}
K_{n,m}
=
\min_{v\in S^{n-1}}
\max_{y\in\Delta_m}
y^\top Gv.
\label{eq:primary}
\end{equation}

We use Gordon's one-sided Gaussian min--max theorem in the
compact-set formulation of Thrampoulidis, Oymak, and Hassibi~\cite{ThrampoulidisOymakHassibi}.
For our choice
$$
\mathcal S_x=S^{n-1},
\qquad
\mathcal S_y=\Delta_m,
\qquad
\psi\equiv0,
$$
it gives
\begin{equation}
\Pp\{K_{n,m}<s\}
\le
2\,\Pp\{\ppphi(g,h)\le s\},
\qquad s\in\R,
\label{eq:gordon}
\end{equation}
where $g\sim N(0,I_m)$ and $h\sim N(0,I_n)$ are independent, and
$$
\ppphi(g,h)
=
\min_{v\in S^{n-1}}
\max_{y\in\Delta_m}
\left\{
g^\top y+\norm{y}\,h^\top v
\right\}.
$$

Since the inner maximum is nondecreasing in the scalar \(h^\top v\), and \[
\min_{v\in S^{n-1}}h^\top v=-\|h\|,
\] the outer minimum is attained at \(v=-h/\|h\|\). 
Therefore
$$
\mathcal A(g,h) = \max_{y\in\Delta_m} \{g^\top y-\|h\|\,\|y\|\}.
 $$

\medskip
\noindent
\textbf{Step 2: scalar characterization of the auxiliary problem.}

Set
$$
r:=\norm{h}.
$$
We consider
\begin{equation}
\mathcal A(g,h)
=
\max_{y\in\Delta_m}
\left\{
g^\top y-r\norm{y}
\right\}.
\label{eq:AO2}
\end{equation}

Let $y^\star$ be an optimizer.  Since $y^\star\in\Delta_m$,
we have $y^\star\neq0$, so the Euclidean norm is differentiable at
$y^\star$.

Introduce a Lagrange multiplier $\nu\in\R$ for the equality constraint
$$
\sum_{i=1}^m y_i=1
$$
and multipliers $\gamma_i\ge0$ for the inequality constraints
$$
y_i\ge0.
$$
Using the Lagrangian
$$
L(y,\nu,\gamma)
=
g^\top y-r\norm{y}
-\nu\left(\sum_{i=1}^m y_i-1\right)
+\sum_{i=1}^m\gamma_i y_i,
$$
the stationarity conditions at $y^\star$ are
\begin{equation}
g_i-\nu+\gamma_i
=
r\,\frac{y_i^\star}{\norm{y^\star}},
\qquad
i=1,\dots,m.
\label{eq:stationarity}
\end{equation}
Complementary slackness gives
\begin{equation}
\gamma_i y_i^\star=0.
\label{eq:slackness}
\end{equation}

If $y_i^\star>0$, then $\gamma_i=0$, and therefore
\begin{equation}
g_i-\nu
=
r\,\frac{y_i^\star}{\norm{y^\star}}.
\label{eq:active}
\end{equation}
In particular, $g_i>\nu$ on every active coordinate.

If $y_i^\star=0$, then \eqref{eq:stationarity} gives $g_i-\nu=-\gamma_i\le0$,
and hence $g_i\le\nu$.
Thus,
$$
y_i^\star>0
\quad\Longleftrightarrow\quad
g_i>\nu.
$$

Consequently,
\begin{equation}
(g_i-\nu)_+
=
r\,\frac{y_i^\star}{\norm{y^\star}},
\qquad i=1,\dots,m.
\label{eq:positive-part}
\end{equation}
Squaring and summing gives
\begin{align}
\sum_{i=1}^m(g_i-\nu)_+^2
&=
\frac{r^2}{\norm{y^\star}^2}
\sum_{i=1}^m(y_i^\star)^2\\
&=
r^2.
\label{eq:root-nu}
\end{align}

It remains to identify $\nu$ with the optimal objective value.
Multiply \eqref{eq:active} by $y_i^\star$ and sum over $i$,
$$
g^\top y^\star
-
\nu\sum_{i=1}^m y_i^\star
=
\frac{r}{\norm{y^\star}}
\sum_{i=1}^m(y_i^\star)^2.
$$
Using $\sum_{i=1}^m y_i^\star=1$ and $\sum_{i=1}^m(y_i^\star)^2=\norm{y^\star}^2$, we obtain
$g^\top y^\star-\nu= r\norm{y^\star}$.
Therefore
$$
\nu
=
g^\top y^\star-r\norm{y^\star}
=
\mathcal A(g,h).
$$
Writing
$
\lambda:=\mathcal A(g,h),
$
equation \eqref{eq:root-nu} becomes
\begin{equation}
\sum_{i=1}^m(g_i-\lambda)_+^2=r^2.
\label{eq:root}
\end{equation}

Thus the optimal value $\lambda$ is the unique solution of
\eqref{eq:root}.  Since the function
 $ s\mapsto
\sum_{i=1}^m(g_i-s)_+^2
$
is continuous and non-increasing, and since $r > 0$ with probability one, we also obtain, for every $s$,
\begin{equation}
\mathcal A(g,h)>s
\quad\Longleftrightarrow\quad
\sum_{i=1}^m(g_i-s)_+^2>\norm{h}^2.
\label{eq:criterion}
\end{equation}

Indeed, equality holds at $s=\lambda$; if $s<\lambda$ the left-hand
side of \eqref{eq:root} is strictly larger than $r^2$, while if
$s>\lambda$ it is strictly smaller.

\medskip
\noindent
\textbf{Step 3: direct integral bounds for the truncated Gaussian moments.}

Let $Z\sim N(0,1)$ and define
$$
I(s):=\E(Z-s)_+^2,
\qquad
J(s):=\E(Z-s)_+^4.
$$
We derive the estimates needed below directly from the Gaussian
integrals, without appealing to Mills-ratio bounds.

By definition,
$$
I(s)
=
\int_s^\infty (z-s)^2\phi(z)\,dz,
$$
where
$$
\phi(z)=\frac{1}{\sqrt{2\pi}}e^{-z^2/2}.
$$
Set $u=z-s$, then
$$I(s) = \int_0^\infty u^2\phi(s+u)\,du
= \phi(s) \int_0^\infty u^2 e^{-su-u^2/2}\,du.
$$
Now set $w=su$. For $s>0$,
\begin{equation}
I(s)
=
\frac{\phi(s)}{s^3}
\int_0^\infty
w^2 e^{-w-w^2/(2s^2)}\,dw.
\label{eq:I-rescaled}
\end{equation}

Moreover,
$$
e^{-w-w^2/(2s^2)}
\le e^{-w},
$$
and therefore
$$
I(s)
\le
\frac{\phi(s)}{s^3}
\int_0^\infty w^2e^{-w}\,dw
=
2\,\frac{\phi(s)}{s^3}. 
$$
For the lower bound, restrict the integral in
\eqref{eq:I-rescaled} to $0\le w\le1$.  Since $s\ge1$,
$$
e^{-w-w^2/(2s^2)}
\ge e^{-3/2}
\qquad
(0\le w\le1),
$$
and hence
$$
I(s)
\ge
\frac{\phi(s)}{s^3}
e^{-3/2}
\int_0^1 w^2\,dw
=
\frac{e^{-3/2}}{3}
\frac{\phi(s)}{s^3}.
$$
Thus, for every $s\ge1$,
\begin{equation}
\frac{e^{-3/2}}{3}\,
\frac{\phi(s)}{s^3}
\le
I(s)
\le
2\,\frac{\phi(s)}{s^3}.
\label{eq:I-bounds}
\end{equation}

Similarly,
$$
J(s)
=
\int_s^\infty (z-s)^4\phi(z)\,dz.
$$
With the same substitutions $u=z-s$ and $w=su$,
$$
J(s)
=
\frac{\phi(s)}{s^5}
\int_0^\infty
w^4 e^{-w-w^2/(2s^2)}\,dw.
$$
For $s\ge1$,
$$
J(s)
\le
\frac{\phi(s)}{s^5}
\int_0^\infty w^4e^{-w}\,dw
=
24\,\frac{\phi(s)}{s^5}.
$$
Hence
\begin{equation}
J(s)
\le
24\,\frac{\phi(s)}{s^5},
\qquad s\ge1.
\label{eq:J-bound}
\end{equation}

We apply these estimates at
$$
t=(1-\rho)\sqrt{2\log\alpha}.
$$
For $\alpha$ sufficiently large, $t\ge1$.  Since
$$
e^{-t^2/2}
=
\alpha^{-(1-\rho)^2},
$$
we have
$$
\phi(t)
=
\frac{1}{\sqrt{2\pi}}
\alpha^{-(1-\rho)^2}.
$$
Therefore, by \eqref{eq:I-bounds},
\begin{align}
\alpha I(t)
&\ge
\frac{e^{-3/2}}{3\sqrt{2\pi}}
\frac{
\alpha^{1-(1-\rho)^2}
}{
t^3
}\\
&=
\frac{e^{-3/2}}{3\sqrt{2\pi}}
\frac{\alpha^\beta}{t^3},
\label{eq:alphaI-lower}
\end{align}
where
$$
\beta
=
1-(1-\rho)^2
=
2\rho-\rho^2>0.
$$
Since
$$
t^3
=
(1-\rho)^3(2\log\alpha)^{3/2},
$$
it follows that
$$
\alpha I(t)\to\infty
\qquad
(\alpha\to\infty).
$$
Hence, for sufficiently large $\alpha$,
\begin{equation}
\alpha I(t)\ge4.
\label{eq:mean-four}
\end{equation}

The same direct estimates also give the ratio needed for the
second-moment bound:
$$
\frac{J(t)}{I(t)^2}
\le
\frac{
24\,\phi(t)t^{-5}
}{
\left(
\frac{e^{-3/2}}{3}\phi(t)t^{-3}
\right)^2
}.
$$
Thus
\begin{equation}
\frac{J(t)}{I(t)^2}
\le
216 e^3\,\frac{t}{\phi(t)}.
\label{eq:J-over-I}
\end{equation}

\medskip
\noindent
\textbf{Step 4: concentration of the truncated Gaussian sum.}

Set
$$
X_i:=(g_i-t)_+^2,
\qquad
S_m:=\sum_{i=1}^mX_i.
$$
Then
$$
\E S_m=mI(t)=n\alpha I(t)\ge4n.
$$

If $S_m\le2n$, then
$$
|S_m-\E S_m|
\ge
\frac12\E S_m.
$$
Chebyshev's inequality along with the independence of the $X_i$ therefore give
\begin{align}
\Pp\{S_m\le2n\}
\le
4\,
\frac{\operatorname{Var}(S_m)}
{(\E S_m)^2}
=
4\,\frac{m(J(t) - I(t)^2)}
{m^2I(t)^2}
\le
4\,
\frac{mJ(t)}
{m^2I(t)^2}.
\label{eq:cheb}
\end{align}

Using \eqref{eq:J-over-I},
$$
\frac{J(t)}{I(t)^2}
\le
216e^3\,\frac{t}{\phi(t)}.
$$
we can estimate
\begin{equation}
\Pp\{S_m\le2n\}
\le
C_1\frac{t}{m\phi(t)}
\label{eq:S-tail}
\end{equation}
for a constant $C_1 >0$.

Since $m=n\alpha$ and
$$
\phi(t)
=
\frac{1}{\sqrt{2\pi}}
\alpha^{-(1-\rho)^2},
$$
we have
\begin{equation}
\frac{t}{m\phi(t)}
=
\sqrt{2\pi}\,
\frac{t}{n\alpha^\beta}
\le
C_2
\frac{\sqrt{\log\alpha}}
{n\alpha^\beta}.
\label{eq:S-explicit}
\end{equation}

On the other hand, since $\norm{h}^2\sim\chi_n^2$, a standard
chi-square upper-tail bound~(See Theorem 2.16 in~\cite{MathDataScience}) gives
\begin{equation}
\Pp\{\norm{h}^2\ge2n\}
\le
e^{-c n}
\label{eq:h-tail}
\end{equation}
for an absolute constant $c>0$.

By \eqref{eq:criterion},
$$
\{\ppphi(g,h)\le t\}
= 
\{S_m\le\norm{h}^2\}
\subseteq
\{S_m\le2n\}
\cup
\{\norm{h}^2\ge2n\}.
$$
Therefore
\begin{equation}
\Pp\{\ppphi(g,h)\le t\}
\le
C_2
\frac{\sqrt{\log\alpha}}
{n\alpha^\beta}
+
e^{-cn}.
\label{eq:AO-tail}
\end{equation}

Gordon's inequality \eqref{eq:gordon} now yields
$$
\Pp\{K_{n,m}<t\}
\le
C
\frac{\sqrt{\log\alpha}}
{n\alpha^\beta}
+
2e^{-cn},
$$
which proves
\eqref{eq:unnormalized-lower-tail}.

\medskip
\noindent
\textbf{Step 5: normalization with an explicit Gaussian norm bound.}

Define
$$
R_{\max}
:=
\max_{1\le i\le m}\norm{g_i}.
$$
For $\rho\in(0,1)$, one-sided Gaussian concentration (see Theorem 13.9 and Proposition 13.15 in \cite{MathDataScience}) gives
$$
\Pp\left\{
\norm{g_i}\ge(1+\rho)\sqrt n
\right\}
\le
e^{-n\rho^2/2}.
$$
A union bound gives
\begin{equation}
\Pp\left\{
R_{\max}\ge(1+\rho)\sqrt n
\right\}
\le
m e^{-n\rho^2/2}.
\label{eq:Rmax-tail}
\end{equation}

On the event
$$
K_{n,m}\ge t>0
$$
we have, for every $v\in S^{n-1}$,
$$
\max_i\ip{g_i}{v}\ge t>0.
$$
Choose $i(v)$ to be the index attaining this maximum.  Then
$$
\max_i\ip{u_i}{v}
\ge
\frac{\max_i\ip{g_i}{v}}
{R_{\max}}.
$$
Taking the minimum over $v$ yields
\begin{equation}
\mu_{n,m}
\ge
\frac{K_{n,m}}{R_{\max}}.
\label{eq:normalization}
\end{equation}

Thus, on the intersection of the events $ K_{n,m}\ge t$ and $R_{\max}\le(1+\rho)\sqrt n$, we obtain
\begin{equation}
\mu_{n,m}
\ge
\frac{1-\rho}{1+\rho}
\sqrt{\frac{2\log\alpha}{n}}.
\label{eq:ratio-bound}
\end{equation}

For $0<\rho<1$,
$$
\frac{1-\rho}{1+\rho}
=
1-\frac{2\rho}{1+\rho}
\ge
1-2\rho.
$$
Since $\rho=\varepsilon/3$,
$$
1-2\rho
=
1-\frac{2\varepsilon}{3}
\ge
1-\varepsilon.
$$
Hence
$$
\mu_{n,m}
\ge
(1-\varepsilon)
\sqrt{\frac{2\log\alpha}{n}}
$$
on this intersection.

Combining
\eqref{eq:unnormalized-lower-tail} and
\eqref{eq:Rmax-tail}, and using $\rho=\frac{\varepsilon}{3}$,
gives
\begin{align}
&\Pp\left\{
\mu_{n,m}
<
(1-\varepsilon)
\sqrt{\frac{2\log\alpha}{n}}
\right\}
\nonumber\\
&\qquad\le
C
\frac{\sqrt{\log\alpha}}
{n\alpha^\beta}
+
2e^{-cn}
+
m\exp\left(-\frac{n\varepsilon^2}{18}\right).
\end{align}
Since $\alpha=m/n$, this is exactly
\eqref{eq:normalized-lower-tail}.

Finally, if
$$
\alpha=\frac mn\to\infty,
\qquad
\log\alpha=o(n),
$$
then
$$
\frac{\sqrt{\log\alpha}}
{n\alpha^\beta}\to0,
\qquad
e^{-cn}\to0,
$$
and
$$
m e^{-n\varepsilon^2/18}
=
\exp\left\{
\log m-\frac{n\varepsilon^2}{18}
\right\}
\to0,
$$
because
$$
\log m
=
\log n+\log\alpha
=
o(n).
$$
The proof is complete.
\end{proof}

\subsubsection{Upper bound for the Gaussian cosine measure}

\begin{theorem}
\label{thm:upper}
There exist absolute constants $C_*>0$ and $\alpha_0>1$ such that,
for every $\alpha\ge\alpha_0$, with
$$
K:=\sqrt{2\log\alpha},
$$
one has
\begin{equation}
\Pp\{K_{n,m}>K\}
\le
1-
\sqrt{
1-\frac{2C_*}{\sqrt{2\pi}\,K}
}.
\label{eq:K-tail}
\end{equation}
Consequently, for every fixed $\varepsilon>0$ and all sufficiently
large $\alpha$,
\begin{align}
&\Pp\left\{
\mu_{n,m}>
(1+\varepsilon)
\sqrt{\frac{2\log(m/n)}{n}}
\right\}
\nonumber\\
&\qquad\le
1-
\sqrt{
1-\frac{2C_*}
{\sqrt{2\pi}\sqrt{2\log(m/n)}}
}
+
m\exp\left\{
-\frac n2
\left(\frac{\varepsilon}{1+\varepsilon}\right)^2
\right\}.
\label{eq:mu-tail}
\end{align}
In particular, if
$$
\frac mn\to\infty,
\qquad
\log\frac mn=o(n),
$$
then, for every fixed $\varepsilon>0$,
$$
\Pp\left\{
\mu_{n,m}>
(1+\varepsilon)
\sqrt{\frac{2\log(m/n)}{n}}
\right\}
\longrightarrow0.
$$
\end{theorem}

The proof combines three standard ingredients: a weighted second-moment construction for perceptron feasibility, following ideas used for the negative spherical perceptron; exponential centering of the one-dimensional weight (these ideas are all inspired by~\cite{MontanariZhongZhou}); and the Hermite-rank/Mehler representation of correlations of functions of Gaussian variables~\cite{BreuerMajor,JansonGaussian}. The specific observation exploited here is that centering annihilates the first Hermite coefficient. Combined with the exact moments of spherical overlaps, this gives a direct second-moment bound without optimizing a nontrivial overlap functional.

\begin{proof}

\medskip
\noindent
\textbf{Step 1: a centered exponential weight.}

Set
$$
K:=\sqrt{2\log\alpha}.
$$
For $c>0$, define
$$
f_c(x):=
e^{-cx}\mathbf 1_{\{x\ge-K\}},
\qquad
p(c):=\E f_c(G),
$$
where $G\sim N(0,1)$.

Writing the expectation as an integral and completing the square,
\begin{align}
p(c)
&=
\frac{1}{\sqrt{2\pi}}
\int_{-K}^{\infty}
e^{-cx-x^2/2}\,dx\\
&=
e^{c^2/2}
\frac{1}{\sqrt{2\pi}}
\int_{-K}^{\infty}
e^{-(x+c)^2/2}\,dx\\
&=
e^{c^2/2}\Phi(K-c).
\label{eq:p}
\end{align}

Noting that for $c \in [0,C]$ the difference quotients of $f_c(G)$ with respect to $c$ are dominated by $e^{CK}|G|$, dominated convergence gives
$$
p'(c)
=
-\E\left[
G e^{-cG}\mathbf 1_{\{G\ge-K\}}
\right],
$$

Thus we have 
$$p'(0) = -\phi(K) < 0,$$
while
$$p(c) \ge e^{cK/2}\Pp \{-K \le G \le -K/2\} \to \infty \\ \text{ as } \\ c \to \infty,$$
guaranteeing the existence of a $c=c_K>0$ with
$$
p'(c)=0.
$$
This is equivalent to
\begin{equation}
\E[Gf_c(G)]=0.
\label{eq:center}
\end{equation}
Differentiating \eqref{eq:p} gives the critical-point equation
\begin{equation}
c\Phi(K-c)=\phi(K-c),
\label{eq:c-equation}
\end{equation}
where
$$
\phi(x)=\frac{1}{\sqrt{2\pi}}e^{-x^2/2}.
$$

Put
$$
t=K-c.
$$
Then
$$
K=t+\frac{\phi(t)}{\Phi(t)}.
$$
Let
$$
f\hspace{-3pt}: \  t\longmapsto
t+\frac{\phi(t)}{\Phi(t)}
$$
and note that $f'(t) = \Var(Z | Z \le t) > 0$ where $Z \sim N(0,1)$. So $f$ is continuous and strictly increasing with $$\lim\limits_{t \to -\infty} f(t) = 0 \text{ and } \lim\limits_{t \to \infty} f(t)= \infty.$$
Hence, for every $K>0$, there is a unique $t$, and therefore a
unique $c=K-t>0$.

As $K\to\infty$, necessarily $t\to\infty$.  From
$$
c=\frac{\phi(t)}{\Phi(t)}
$$
we obtain
$$
c\to0
$$
and
$$
tc\le 2t\phi(t)\to0
$$
for all sufficiently large $t$, since $\Phi(t)\ge1/2$.
Consequently
\begin{equation}
Kc=tc+c^2\longrightarrow0.
\label{eq:Kc}
\end{equation}

It will be convenient in our quantitative estimates below to use the following bound: There exists
$K_0$ such that, for all $K\ge K_0$,
\begin{equation}
Kc\le\frac{1}{8e^2}.
\label{eq:Kc-explicit}
\end{equation}

Define
$$
f:=f_{c_K},
\qquad
p:=\E f(G),
\qquad
F:=\frac{f}{p},
\qquad
V:=\operatorname{Var}(F(G)).
$$
Then
\begin{equation}
\E F(G)=1,
\qquad
\E[GF(G)]=0,
\label{eq:F-center}
\end{equation}
where the second equality follows from \eqref{eq:center}.

\medskip
\noindent
\textbf{Step 2: an explicit variance bound from direct Gaussian
integrals.}

Let
$$
r:=\Phi(-K)
=
\int_K^\infty\phi(x)\,dx.
$$
With $x=K+u$,
$$
r
=
\phi(K)
\int_0^\infty
e^{-Ku-u^2/2}\,du.
$$
Setting $w=Ku$ gives
\begin{equation}
r
=
\frac{\phi(K)}{K}
\int_0^\infty
e^{-w-w^2/(2K^2)}\,dw.
\label{eq:r-integral}
\end{equation}
For $K\ge1$, restricting to $0\le w\le1$ yields
\begin{equation}
r
\ge
e^{-3/2}\frac{\phi(K)}{K}.
\label{eq:r-lower}
\end{equation}

For $a=1,2$, define the shifted tails
$$
r_a
:=
1-\Phi(K-ac)
=
\int_{K-ac}^{\infty}\phi(x)\,dx.
$$
Since $K-ac<K$,
$$
r_a-r
=
\int_{K-ac}^{K}\phi(x)\,dx.
$$
For sufficiently large $K$, the Gaussian density is decreasing
throughout this interval, so
\begin{equation}
0\le r_a-r
\le
ac\,\phi(K-ac).
\label{eq:shift-difference}
\end{equation}

Furthermore,
\begin{equation}
\frac{\phi(K-ac)}{\phi(K)}
=
\exp\left(
aKc-\frac{a^2c^2}{2}
\right)
\le e^{aKc}.
\label{eq:phi-ratio}
\end{equation}
Using \eqref{eq:Kc-explicit}, for $a=1,2$,
$$
\phi(K-ac)\le e^{1/(4e^2)}\phi(K)<e^{1/2}\phi(K).
$$
Combining this with \eqref{eq:r-lower} and using \eqref{eq:Kc-explicit} gives
\begin{align*}
\frac{r_a-r}{r}
&\le
a c e^{1/2}\phi(K)
\frac{K e^{3/2}}{\phi(K)}\\
&=
a e^2 Kc\\
&\le\frac a8.
\end{align*}
Hence by \eqref{eq:shift-difference} and the above we obtain
\begin{equation}
r\le r_1\le\frac98 r,
\qquad
r\le r_2\le\frac54 r.
\label{eq:shifted-tail-bounds}
\end{equation}

We also need an explicit bound on $c^2$ relative to $r$.
From \eqref{eq:c-equation},
$$
c=\frac{\phi(K-c)}{\Phi(K-c)}.
$$
For large $K$, $K-c>0$, and hence $\Phi(K-c)\ge1/2$.
So by \eqref{eq:phi-ratio}
$$
c
\le
2\phi(K-c)
\le
2e^{Kc}\phi(K).
$$
Using \eqref{eq:Kc-explicit},
$$
c\le3\phi(K)
$$
for all sufficiently large $K$.  Thus, by \eqref{eq:r-lower},
$$
\frac{c^2}{r}
\le
9e^{3/2}K\phi(K).
$$
Since $K\phi(K)\to0$, after increasing $K_0$ if necessary we may
assume
\begin{equation}
c^2\le\frac18 r.
\label{eq:c2}
\end{equation}

Now
$$
p=\E f(G)=e^{c^2/2}\Phi(K-c).
$$
Recall that
\[
f(x)=e^{-cx}\mathbf 1_{\{x\ge -K\}},
\]
so
\[
f(x)^2=e^{-2cx}\mathbf 1_{\{x\ge -K\}},
\]
Thus
\[
\begin{aligned}
\mathbb E[f(G)^2]
&=\frac1{\sqrt{2\pi}}\int_{-K}^{\infty}
e^{-2cx-x^2/2}\,dx.\\
&=e^{2c^2}\frac1{\sqrt{2\pi}}
\int_{-K+2c}^{\infty}e^{-u^2/2}\,du\\
&=e^{2c^2}\bigl[1-\Phi(-K+2c)\bigr]\\
&=e^{2c^2}\Phi(K-2c).
\end{aligned}
\]
Noting that  \[
V=\mathbb E[F(G)^2]-\bigl(\mathbb E[F(G)]\bigr)^2
=\frac{\mathbb E[f(G)^2]}{p^2} -1,
\]
we have
\begin{equation}
1+V
=
e^{c^2}
\frac{\Phi(K-2c)}{\Phi(K-c)^2}
=
e^{c^2}
\frac{1-r_2}{(1-r_1)^2}.
\label{eq:V-exact}
\end{equation}

For large $K$ we may also assume $r\le1/8$.  Then
$$
r_1\le\frac98r<\frac14.
$$
Also for $0\le x\le1/4$,
$$
(1-x)^{-2}\le1+4x.
$$
Moreover, since $c^2\le r/8$,
$$
e^{c^2}\le1+2c^2\le1+\frac r4.
$$
Thus, using $1-r_2\le1$ in \eqref{eq:V-exact},
\begin{align*}
1+V
&\le
\left(1+\frac r4\right)
(1+4r_1)\\
&\le
\left(1+\frac r4\right)
\left(1+\frac92r\right).
\end{align*}
For $r\le1/8$, the right-hand side is at most $1+6r$.
Therefore
\begin{equation}
V\le C_*\,\Phi(-K).
\label{eq:V-bound}
\end{equation}
For example, we can choose $C_*=6$.

From \eqref{eq:r-integral}
$$
\Phi(-K)\le\frac{\phi(K)}{K}.
$$
Since $K=\sqrt{2\log\alpha}$ and  $e^{-K^2/2}=\alpha^{-1}$, we obtain
\begin{equation}
\alpha V
\le
\frac{C_*}{\sqrt{2\pi}\,K}.
\label{eq:alphaV}
\end{equation}

\medskip
\noindent
\textbf{Step 3: centering gives a global quadratic correlation bound.}

Let $(G_1,G_2)$ be jointly standard Gaussian with correlation
$q\in(-1,1)$.  We use the Hermite polynomials in  the ``probabilist's normalization'', namely
$$
H_j(x)
=
(-1)^j e^{x^2/2}
\frac{d^j}{dx^j}e^{-x^2/2},
$$
which satisfy
$$
\E[H_j(G)H_k(G)]
=
j!\,\mathbf 1_{\{j=k\}}.
$$

Expand
$$
F(x)
=
\sum_{j=0}^{\infty}
\frac{a_j}{j!}H_j(x).
$$
By \eqref{eq:F-center}, we have $a_0=1$ and $a_1=0$.

Recall Mehler's identity (\cite{NIST:DLMF} 18.18.28) for this normalization (see 18.7.11 there)

\[
\frac{1}{\sqrt{1-q^2}}
\exp\!\left(
\frac{qxy-\frac12q^2(x^2+y^2)}{1-q^2}
\right)
=
\sum_{j=0}^{\infty}\frac{q^j}{j!}H_j(x)H_j(y).
\]
Multiplying each side by the product of two independent standard Gaussian pdfs gives the following formula for the joint distribution
\[
\phi_q(x,y)
=
\phi(x)\phi(y)
\sum_{j=0}^{\infty}\frac{q^j}{j!}H_j(x)H_j(y).
\]
This gives us the following identity, since $F$ is square-integrable with respect to the standard Gaussian measure:
\begin{align*}
\E[F(G_1)F(G_2)]
&=
\int_{\mathbb R^2}F(x)F(y)\phi_q(x,y)\,dx\,dy \\
&=
\sum_{j=0}^{\infty}\frac{q^j}{j!}
\left(\int_{\mathbb R}F(x)H_j(x)\phi(x)\,dx\right)
\left(\int_{\mathbb R}F(y)H_j(y)\phi(y)\,dy\right) \\
&=
1+
\sum_{j=2}^{\infty}
\frac{a_j^2}{j!}q^j
\end{align*}

Since, for $j\ge2$ and $|q|\le1$,
$$
q^j\le |q|^j\le q^2,
$$
we obtain
\begin{align}
\E[F(G_1)F(G_2)]
&\le
1+
q^2
\sum_{j=2}^{\infty}
\frac{a_j^2}{j!}.
\end{align}
By Parseval's identity,
$$
\sum_{j=2}^{\infty}
\frac{a_j^2}{j!}
=
\operatorname{Var}(F(G))
=
V.
$$
Hence
\begin{equation}
\E[F(G_1)F(G_2)]
\le
1+Vq^2.
\label{eq:global-correlation}
\end{equation}
The same bound holds at \(q=\pm1\) by Cauchy–Schwarz.

\medskip
\noindent
\textbf{Step 4: second moment of the weighted feasible volume.}

Let $\sigma$ denote normalized surface measure on $S^{n-1}$ and set
$$
Z
:=
\int_{S^{n-1}}
\prod_{i=1}^m
f(\ip{g_i}{v})\,d\sigma(v).
$$
Since $f(x)=0$ when $x<-K$, the event $Z>0$ implies the existence of
$v\in S^{n-1}$ such that
\begin{equation}
\ip{g_i}{v}\ge-K,
\qquad i=1,\dots,m.
\label{eq:feasible}
\end{equation}

For each fixed $v$, the random variables
$$
\ip{g_1}{v},\dots,\ip{g_m}{v}
$$
are independent standard Gaussians, as they are linear combinations of such Gaussians.  Hence
$$
\E Z=p^m.
$$

Let $V_1,V_2$ be independent uniform random vectors on $S^{n-1}$ and
set
$$
Q:=\ip{V_1}{V_2}.
$$
Conditioned on $Q=q$, each pair
$$
\bigl(
\ip{g_i}{V_1},
\ip{g_i}{V_2}
\bigr)
$$
is jointly standard Gaussian with correlation $q$.  Therefore
\begin{align*}
\frac{\E Z^2}{(\E Z)^2}
&=
\frac{\E_Q
\left[
\E[f(G_1)f(G_2)]^m
\right]}{p^{2m}}\\
&=
\E_Q
\left[
\E[F(G_1)F(G_2)]^m
\right]\\
&\le
\E_Q(1+VQ^2)^m\\
&\le
\E_Qe^{mVQ^2},
\label{eq:second-moment}
\end{align*}
where in the last two lines we used \eqref{eq:global-correlation} and the inequality $1+x\le e^x$.

We need the following lemma in the next calculation. Its proof can be found in \Cref{sec:appendix}.

\begin{restatable}{lemma}{evenmoments}
\label{lem:even-moments}
The even moments of $Q$ are
$$
\E Q^{2j}
=
\frac{(2j-1)!!}
{n(n+2)\cdots(n+2j-2)}
\le
\frac{(2j-1)!!}{n^j}.
$$

\end{restatable}

By this lemma, whenever $2mV<n$,
\begin{align*}
\E e^{mVQ^2}
&\le
\sum_{j=0}^{\infty}
\frac{(mV)^j}{j!}
\frac{(2j-1)!!}{n^j}\\
&=
\left(
1-\frac{2mV}{n}
\right)^{-1/2}\\
&=
(1-2\alpha V)^{-1/2},
\end{align*}
where we used the binomial series expansion for $(1-x)^{-1/2}$.

Using \eqref{eq:alphaV}, for sufficiently large $\alpha$ we have
$$
2\alpha V<1,
$$
and
\begin{equation}
\frac{\E Z^2}{(\E Z)^2}
\le
\left(
1-\frac{2C_*}{\sqrt{2\pi}\,K}
\right)^{-1/2}.
\label{eq:second-explicit}
\end{equation}

The Paley--Zygmund inequality gives
$$
\Pp(Z>0)
\ge
\frac{(\E Z)^2}{\E Z^2},
$$
and therefore
\begin{equation}
\Pp(Z>0)
\ge
\sqrt{
1-\frac{2C_*}{\sqrt{2\pi}\,K}
}.
\label{eq:Z-positive}
\end{equation}

On the event $Z>0$, choose $v$ satisfying \eqref{eq:feasible}.
Replacing $v$ by $-v$ gives
$$
\ip{g_i}{-v}\le K
\qquad
\text{for all }i.
$$
Hence
$$
K_{n,m}\le K.
$$
Consequently
$$
\Pp\{K_{n,m}>K\}
\le
1-
\sqrt{
1-\frac{2C_*}{\sqrt{2\pi}\,K}
},
$$
which proves \eqref{eq:K-tail}.

\medskip
\noindent
\textbf{Step 5: normalization with an explicit Gaussian norm tail.}

Fix $\eta\in(0,1)$ and define
$$
R_{\min}:=
\min_{1\le i\le m}\norm{g_i}.
$$
The standard Gaussian norm lower-tail estimate gives
$$
\Pp\left\{
\norm{g_i}\le(1-\eta)\sqrt n
\right\}
\le
e^{-n\eta^2/2}.
$$
By the union bound,
\begin{equation}
\Pp\left\{
R_{\min}\le(1-\eta)\sqrt n
\right\}
\le
m e^{-n\eta^2/2}.
\label{eq:Rmin}
\end{equation}

On the event
$$
K_{n,m}\le K,
\qquad
R_{\min}>(1-\eta)\sqrt n,
$$
there exists $v\in S^{n-1}$ such that
$$
\ip{g_i}{v}\le K
\qquad
\text{for every }i.
$$
Then
$$
\ip{u_i}{v}
=
\frac{\ip{g_i}{v}}{\norm{g_i}}
\le
\frac{K}{\norm{g_i}}
\le
\frac{K}{(1-\eta)\sqrt n}.
$$

Hence
$$
\mu_{n,m}
\le
\frac{1}{1-\eta}
\sqrt{\frac{2\log\alpha}{n}}.
$$

Combining \eqref{eq:K-tail} and \eqref{eq:Rmin} and union bounding gives,
\begin{align*}
&\Pp\left\{
\mu_{n,m}>
\frac{1}{1-\eta}
\sqrt{\frac{2\log\alpha}{n}}
\right\}
\nonumber\\
&\qquad\le
\mathbb P\{K_{n,m}>K\}
+
\mathbb P\{R_{\min}\le(1-\eta)\sqrt n\}\\
&\qquad\le
1-
\sqrt{
1-\frac{2C_*}
{\sqrt{2\pi}\sqrt{2\log\alpha}}
}
+
m e^{-n\eta^2/2}.
\label{eq:eta-tail}
\end{align*}

Given $\varepsilon>0$, choose
$$
\eta=\frac{\varepsilon}{1+\varepsilon}.
$$
Then
$$
\frac{1}{1-\eta}=1+\varepsilon,
$$
which gives \eqref{eq:mu-tail}.

Finally, if
$$
\alpha=\frac mn\to\infty,
\qquad
\log\alpha=o(n),
$$
then the first term on the right-hand side of \eqref{eq:mu-tail}
tends to zero because
$
\sqrt{\log\alpha}\to\infty$,
while
$$
m e^{-n\eta^2/2}
=
\exp\left\{
\log m-\frac{n\eta^2}{2}
\right\}
\to0
$$
because
$
\log m
=
\log n+\log\alpha
=
o(n)$.

\end{proof}

\end{section}

\section*{Acknowledgements}

T.S.\ wants to thank Benjamin Friedlander, whose inquisitive questions about conic combinations inspired the research presented in this paper, and March Boedihardjo for helpful discussions on conic combinations.
The authors acknowledge support from NIH R01HL163515, NSF DMS-2208356, and DOE DE-SC0023490.

\bibliography{conic}

\appendix

\section{Proof of Lemma 5.5}
\label{sec:appendix}

\evenmoments*

\begin{proof}
    First, reduce the inner product to one coordinate. Condition on \(V_1=v\), and choose an orthogonal matrix \(O\) sending \(v\) to \(e_1\). Then
\[
\langle v,V_2\rangle
=\langle Ov,OV_2\rangle
=\langle e_1,OV_2\rangle.
\]
Since \(OV_2\) is still uniform on the sphere, this has the same distribution as the first coordinate \(U_1\) of a uniform vector \(U\). This holds for every \(v\), so
\[
Q\overset{d}=U_1.
\]
Next, express \(U\) using a Gaussian vector. Let
\[
X=(X_1,\ldots,X_n)\sim N(0,I_n),
\qquad R=\|X\|,
\qquad U=\frac XR.
\]
Then \(U\) is uniform on the sphere and independent of \(R\). These facts follow from the Gaussian density in polar coordinates: with respect to \(dr\,d\sigma(u)\), it is
\[
C_n r^{n-1}e^{-r^2/2},
\]
which factors into a density depending only on \(r\) and uniform surface measure in \(u\).
Because \(X_1=RU_1\), independence gives
\[
\mathbb E X_1^{2j}
=\mathbb E[R^{2j}U_1^{2j}]
=\mathbb E R^{2j}\,\mathbb E U_1^{2j}.
\]
Consequently,
\begin{equation}
\mathbb E Q^{2j}
=\mathbb E U_1^{2j}
=\frac{\mathbb E X_1^{2j}}{\mathbb E R^{2j}}.
\label{eqn:EQ}
\end{equation}
We can calculate both moments by integration by parts.
For the numerator, use \(\phi'(x)=-x\phi(x)\):
\[
\begin{aligned}
\mathbb E X_1^{2j}
&=\int_{\mathbb R}x^{2j}\phi(x)\,dx\\
&=-\int_{\mathbb R}x^{2j-1}\phi'(x)\,dx\\
&=(2j-1)\int_{\mathbb R}x^{2j-2}\phi(x)\,dx\\
&=(2j-1)\mathbb E X_1^{2j-2}.
\end{aligned}
\]
The boundary term vanishes because the Gaussian density decays faster than any polynomial. Starting from \(\mathbb E X_1^0=1\), this recurrence gives
\[
\mathbb E X_1^{2j}=(2j-1)!!.
\]
For the denominator, the radial density above gives
\[
\mathbb E R^{2j}
=C_n\int_0^\infty r^{n+2j-1}e^{-r^2/2}\,dr.
\]
Write the integrand as \(r^{n+2j-2}(r e^{-r^2/2})\), and use
\[
\frac{d}{dr}e^{-r^2/2}=-r e^{-r^2/2}.
\]
Integration by parts yields
\[
\begin{aligned}
\mathbb E R^{2j}
&=(n+2j-2)C_n
  \int_0^\infty r^{n+2j-3}e^{-r^2/2}\,dr\\
&=(n+2j-2)\mathbb E R^{2j-2}.
\end{aligned}
\]
Again the boundary terms vanish. Starting from \(\mathbb E R^0=1\), we obtain
\[
\mathbb E R^{2j}
=n(n+2)\cdots(n+2j-2).
\]
Substituting into \eqref{eqn:EQ} proves the formula.
Finally, every denominator factor is at least \(n\), so
\[
\mathbb E Q^{2j}\le\frac{(2j-1)!!}{n^j}.
\]
For \(j=0\), both sides are \(1\), using the conventions that an empty product is \(1\) and \((-1)!!=1\).
\end{proof}

\section{Explanation of the inaccuracy in Bonnet and O'Reilly}
\label{sec:explain-inacc}
Using the notation in~\cite[Lemma~6.10]{BonnetOReilly2022}, with \(n\) points in dimension \(d\), Bonnet and O'Reilly claim there exist constants $C_1, C_2$ such that:
\[
C_1\frac{r_1d\ln(n/d)}{n}
\;\boldsymbol{\ge}\;
C_2\frac{r_1d\ln n}{n-d}.
\]

We now demonstrate that for a certain choice of $n$ this does not hold. Choose
\[
n=\lceil d\log d\rceil,\qquad d\to\infty.
\]
The ceiling operation introduces an error smaller than \(1\), so we can write
\[
n=d\log d+\delta_d,\qquad 0\le\delta_d<1.
\]
Consequently,
\[
\frac nd=\log d+\frac{\delta_d}{d}.
\]
The basic limit used below is
\[
\frac{\log\log d}{\log d}\longrightarrow0.
\]
Indeed, setting \(x=\log d\), this becomes the familiar limit \(\log x/x\to0\).
1. The first expression satisfies
   \[
   \frac{d\log(n/d)}n
   \sim\frac{\log\log d}{\log d}
   \longrightarrow0.
   \]
   First divide its numerator and denominator by \(d\):
   \[
   \frac{d\log(n/d)}n
   =
   \frac{\log(n/d)}{n/d}.
   \]
   We already know that
   \[
   \frac nd=\log d+\frac{\delta_d}{d}\sim\log d.
   \]
   For the numerator, factor \(\log d\) inside the logarithm:
   \[
   \begin{aligned}
   \log(n/d)
   &=\log\left(\log d+\frac{\delta_d}{d}\right)\\
   &=\log\left[
      \log d\left(1+\frac{\delta_d}{d\log d}\right)
      \right]\\
   &=\log\log d+
     \log\left(1+\frac{\delta_d}{d\log d}\right).
   \end{aligned}
   \]
   The final term tends to zero because
   \[
   0\le\frac{\delta_d}{d\log d}
   <\frac1{d\log d}\longrightarrow0.
   \]
   Thus
   \[
   \log(n/d)=\log\log d+o(1)\sim\log\log d.
   \]
   Combining the numerator and denominator estimates gives
   \[
   \frac{\log(n/d)}{n/d}
   \sim
   \frac{\log\log d}{\log d}
   \longrightarrow0.
   \]
   Here \(f_d\sim g_d\) means \(f_d/g_d\to1\).
2. The second expression satisfies
   \[
   \frac{d\log n}{n-d}\longrightarrow1.
   \]
   Again divide the numerator and denominator by \(d\):
   \[
   \frac{d\log n}{n-d}
   =
   \frac{\log n}{n/d-1}.
   \]
   For the numerator, use \(n=d(n/d)\):
   \[
   \begin{aligned}
   \log n
   &=\log d+\log(n/d)\\
   &=\log d+\log\log d+o(1).
   \end{aligned}
   \]
   The denominator is
   \[
   \frac nd-1
   =\log d-1+\frac{\delta_d}{d}
   \]
   Therefore,
   \[
   \frac{d\log n}{n-d}
   =
   \frac{\log d+\log\log d+o(1)}
        {\log d-1+\delta_d/d}.
    =  
   \frac{
      1+\dfrac{\log\log d}{\log d}
       +\dfrac{o(1)}{\log d}
   }{
      1-\dfrac1{\log d}
       +\dfrac{\delta_d}{d\log d}
   } \to 1.
   \]
   Thus the desired constants cannot exist.
\end{document}